\documentclass[a4paper]{amsart}
\title{The geometry of Fronts}
\date{December 09, 2006}
\usepackage[dvips]{graphicx} 
\usepackage{verbatim,enumerate}
\usepackage{amssymb}
\usepackage[usenames]{color}
\usepackage{amsthm}
\theoremstyle{plain}
 \newtheorem{theorem}{Theorem}[section]
 \newtheorem*{theorem*}{Theorem}
 \newtheorem*{lemma*}{Lemma}
 \newtheorem{proposition}[theorem]{Proposition}
 \newtheorem{fact}[theorem]{Fact}
 \newtheorem{fact*}{Fact}
 \newtheorem{lemma}[theorem]{Lemma}
 \newtheorem{corollary}[theorem]{Corollary}
\theoremstyle{remark}
 \newtheorem{definition}[theorem]{Definition}
 \newtheorem{remark}[theorem]{Remark}
 \newtheorem*{remark*}{Remark}
 \newtheorem*{acknowledgements}{Acknowledgements}
 \newtheorem{example}[theorem]{Example}
\numberwithin{equation}{section}

\newcommand{\Z}{\boldsymbol{Z}}
\newcommand{\R}{\boldsymbol{R}}
\newcommand{\C}{\boldsymbol{C}}

\newcommand{\sech}{\operatorname{sech}}
\renewcommand{\phi}{\varphi}

\newcommand{\sign}{\operatorname{sgn}}
\newcommand{\ext}{\operatorname{ext}}
\newcommand{\inner}[2]{\left\langle{#1},{#2}\right\rangle}
\newcommand{\Area}{\operatorname{Area}}

\newcommand{\E}{\mathcal{E}}
\newcommand{\PSL}{\operatorname{PSL}}

\author{Kentaro Saji}
\address[Saji]{%
   Department of Mathematics,
   Hokkaido University,
   Sapporo 060-0810,
   Japan
}
\email{saji@math.sci.hokudai.ac.jp}

\author{Masaaki Umehara}
\address[Umehara]{%
   Department of Mathematics, Graduate School of Science,
   Osaka University,
   Toyonaka, Osaka 560-0043,
   Japan
}
\email{umehara@math.wani.osaka-u.ac.jp}
\author{Kotaro Yamada}
\address[Yamada]{%
   Faculty of Mathematics,
   Kyushu University,
   Higashi-ku, Fukuoka 812-8581, Japan%
}
\email{kotaro@math.kyushu-u.ac.jp}
\begin{document}
\begin{abstract}
 We shall introduce the {\em singular curvature function\/}
 on cuspidal edges of surfaces,
 which is related to the Gauss-Bonnet formula 
 and which
 characterizes the shape of cuspidal edges.
 Moreover, it is closely related to the behavior of 
 the Gaussian curvature of a surface near cuspidal edges
 and swallowtails.
\end{abstract}
\maketitle
\section*{Introduction}
\begin{figure}[h]
\renewcommand{\thefigure}{\arabic{figure}}
 \begin{center}
   \begin{tabular}{c@{\hspace{1.5cm}}c}
        \includegraphics[width=3.5cm]{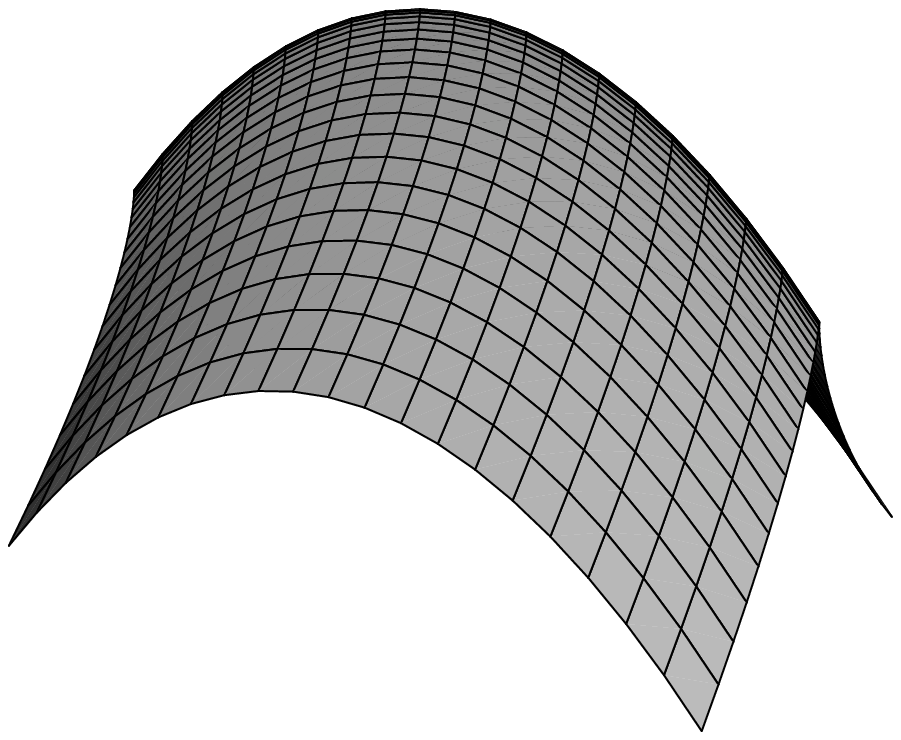}&
        \includegraphics[width=3.5cm]{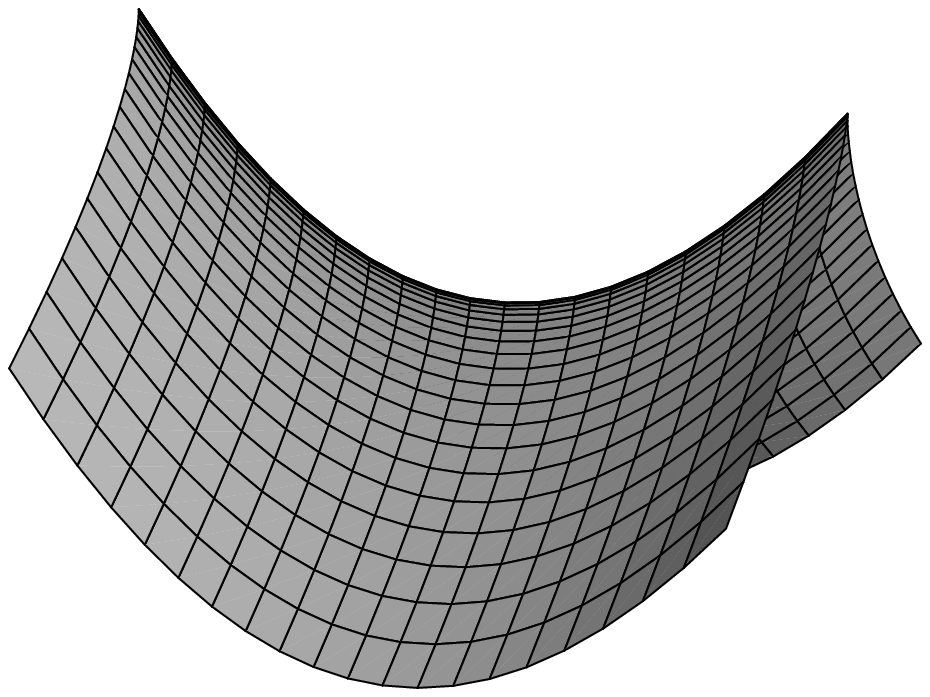}
   \end{tabular}
 \end{center}
\caption{Positively and negatively curved cuspidal edges 
           (Example~\ref{ex:parabola}).}
\label{fig:parabola}
\end{figure}

Let $M^2$ be an oriented $2$-manifold and $f\colon{}M^2\to \R^3$ 
a $C^\infty$-map. 
A point $p\in M^2$ is called a {\em singular point\/} 
if $f$ is not an immersion at $p$.
A singular point is called a 
{\em cuspidal edge\/} or {\em swallowtail\/} 
if it is locally diffeomorphic to 
\begin{equation}\tag{1}\label{eq:cuspidal-swallow}
 f_C(u,v):=(u^2,u^3,v) \quad\text{or}\quad
 f_S(u,v):=(3u^4+u^2v,4u^3+2uv,v)
\end{equation}
at $(u,v)=(0,0)$, respectively.
These two types of singular points characterize 
the generic singularities of wave fronts 
(cf.\ \cite{AGV};
 for example, parallel surfaces of immersed surfaces in $\R^3$ are
 fronts), 
and we have a useful criterion 
(Fact~\ref{fact:intrinsic-criterion}; cf.\ \cite{KRSUY}) 
for determining them.
It is of interest to investigate these singularities
from the viewpoint of differential geometry.
In this paper, we shall distinguish two types of cuspidal edges
as in Figure~\ref{fig:parabola}.
More precisely, we shall define the singular curvature function
$\kappa_s$ along cuspidal edges. 
The left-hand figure in Figure~\ref{fig:parabola} is 
positively curved and the right-hand figure
is negatively curved
(see Corollary~\ref{cor:positive-negative-cuspidal-edge}).

The definition of the singular curvature function does not depend 
on the orientation nor on the co-orientation of the front and is 
closely related to the following two Gauss-Bonnet formulas 
given by Langevin-Levitt-Rosenberg and Kossowski 
when $M^2$ is compact:
\begin{alignat*}{2}
 2\deg(\nu)&= \chi(M_+)-\chi(M_-)
         +\#S_+-\#S_- \qquad &&(\text{\cite{LLR},\cite{K1}})
 \tag{2}\label{eq:GB-signed}\\
 2\pi\chi(M^2)
         &=\int_{M^2}K\,dA+2\int_{\text{Singular set}}\kappa_s\, ds
               \qquad &&(\text{\cite{K1}}),
 \tag{3}\label{eq:GB-unsigned}
\end{alignat*}
where $\deg(\nu)$ is the degree of the Gauss map $\nu$,
$\#S_+,\#S_-$ are the numbers of positive and negative 
swallowtails respectively (see Section~\ref{sec:GB}), 
and $M_+$ (resp.\ $M_-$) is the open submanifold of $M^2$
to which the co-orientation is compatible 
(resp.\ not compatible)
with respect to the orientation.
In the proofs of these formulas in \cite{LLR} and \cite{K1}, 
the singular curvature implicitly appeared as 
a form  $\kappa_s\,ds$.
(Formula \eqref{eq:GB-signed} stated in \cite{LLR}, and 
  proofs for both \eqref{eq:GB-signed} and \eqref{eq:GB-unsigned}
  are in \cite{K1}.)

Recently, global properties of fronts were investigated 
via flat surfaces in hyperbolic $3$-space $H^3$ (\cite{KUY1,KRSUY}), 
via maximal surfaces in Minkowski $3$-space (\cite{UY}), and 
via constant mean curvature one surfaces in de Sitter space 
(\cite{F}, see also Lee and Yang \cite{LY}).
Such surfaces satisfy certain Osserman type inequalities 
for which equality characterizes the proper embeddedness of
their ends.
We also note that Mart\'\i{}nez \cite{Mar} investigated
global properties of improper affine spheres with singularities, 
which are related to flat fronts in $H^3$. 
(See also Ishikawa and Machida \cite{IM}.)

The purpose of this paper is to give geometric meaning to 
the singular curvature function and investigate its properties. 
For example, it diverges to $-\infty$ at swallowtails
(Corollary~\ref{cor:singular-curvature-peak}).
Moreover, we shall investigate behavior 
of the Gaussian curvature $K$ near singular points. 
For example, the Gaussian curvature $K$ is generically unbounded near 
cuspidal edges and swallowtails and will take different signs
from the left-hand side to the right-hand side of 
a singular curve. 
However, on the special occasions that $K$ is bounded,
the shape of these singularities is very restricted:
for example, 
{\it singular curvature is non-positive if the Gaussian curvature
     is non-negative} 
(Theorem~\ref{thm:gaussian-singular}).
A similar phenomena holds for the case of hypersurfaces 
(Section~\ref{sec:hyper}).

The paper is  organized as follows:
In Section~\ref{sec:curvature}, we define the singular curvature,
and give its fundamental properties.
In Section~\ref{sec:GB}, we generalize the two Gauss-Bonnet formulas
\eqref{eq:GB-signed} and \eqref{eq:GB-unsigned} 
to fronts which admit finitely many 
corank one ``peak'' singularities.
In Section~\ref{sec:gaussian-curvature}, 
we investigate behavior of Gaussian curvature.
Section~\ref{sec:zigzag} is devoted to formulating a 
topological invariant of closed fronts called  the
``zig-zag number''
(introduced in \cite{LLR}) 
from the viewpoint of differential geometry.
We shall generalize the results of 
Section~\ref{sec:gaussian-curvature} to hypersurfaces
in Section~\ref{sec:hyper}.
Finally, in Section~\ref{sec:intrinsic}, we introduce an intrinsic
formulation of the geometry of fronts.
\begin{acknowledgements}
The authors thank Shyuichi Izumiya, Go-o Ishikawa,  
Osamu Saeki, Osamu Kobayashi
and Wayne Rossman for fruitful discussions and valuable comments.
\end{acknowledgements}
\section{Singular curvature}
\label{sec:curvature}
Let $M^2$ be an oriented $2$-manifold and $(N^3,g)$ an oriented
Riemannian $3$-manifold.
The unit cotangent bundle $T_1^*N^3$ has the canonical contact structure
and can be identified with the unit tangent bundle $T_1 N^3$.
A smooth map $f\colon{}M^2 \to N^3$ is called a {\em front\/}
if there exists a unit vector field $\nu$ of $N^3$ along $f$
such that $L:=(f,\nu)\colon{}M^2\to T_1N^3$ is a Legendrian immersion
(which is also called an isotropic immersion),
that is, the pull-back of the canonical contact form of $T_1N^3$
vanishes on $M^2$.
This condition is equivalent to the following orthogonality condition:
\begin{equation}\label{eq:orthogonal}
   g(f_*X,\nu)=0 \qquad (X\in TM^2),
\end{equation}
where $f_*$ is the differential map of $f$.
The vector field  $\nu$ is called the {\it unit normal vector\/} of
the front $f$.
The {\em first fundamental form $ds^2$}
and the {\em second fundamental form $h$} of the front are defined 
in the same way as for surfaces:
\begin{equation}\label{eq:first-second}
 ds^2(X,Y) := g(f_*X,f_*Y),~
 h(X,Y)    := -g(f_*X,D_Y\nu)\qquad
 \bigl(X,Y\in TM^2\bigr),
\end{equation}
where $D$ is the Levi-Civita connection of $(N^3,g)$.

We denote by $\mu_g$ the Riemannian volume element of $(N^3,g)$.
Let $f\colon{}M^2\to N^3$ be a front and $\nu$ the unit normal vector of
$f$, and set
\begin{equation}\label{eq:signed-area-form}
  d\hat A:=f^*(\iota_\nu\mu_g)=
        \mu_g(f_u,f_v,\nu)\,du\wedge dv
	\quad 
	\left(f_u=f_*\left(\frac{\partial}{\partial u}\right),
              f_v=f_*\left(\frac{\partial}{\partial v}\right)
        \right),	      
\end{equation}
called the {\it signed area form\/}, where 
$(u,v)$ is a local coordinate system of $M^2$ and 
$\iota_\nu$ is the interior product with respect to $\nu\in TN^3$.
Suppose now that $(u,v)$ is compatible to the orientation of $M^2$.
Then the function
\begin{equation}\label{eq:signed-area-density}
   \lambda(u,v):=\mu_g(f_u,f_v,\nu)
\end{equation}
is called the ({\em local}) {\em signed area density function}. %
We also  set 
\begin{multline}\label{eq:absolute-area-form}
 d A:=|\mu_g(f_u,f_v,\nu)|\,du\wedge dv
        =\sqrt{EG-F^2}\,du\wedge dv = |\lambda|\,du\wedge\,dv\\
   \bigl(E:=g(f_u,f_u), F:=g(f_u,f_v), G:=g(f_v,f_v)\bigr),
\end{multline}
which is independent of the choice of orientation-compatible coordinate
system $(u,v)$ and is called the 
({\em absolute}) {\em area form\/} of $f$.
Let $M_+$ (resp.\ $M_-$) be the open submanifolds where 
the ratio $(d\hat A)/(dA)$ is positive (resp.\ negative). 
If $(u,v)$ is a coordinate system compatible to the orientation of
$M^2$, the point $(u,v)$ belongs to $M_+$ (resp.\ $M_-$)
if and only if $\lambda(u,v)>0$  ($\lambda(u,v)<0$), 
where $\lambda$ is the signed area density function.
\begin{definition}\label{def:nondeg}
 Let $f\colon{}M^2\to N^3$ be a front. 
 A point $p\in M^2$ is called a {\it singular point\/} if $f$ is not an
 immersion at $p$.
 We call the set of singular points of $f$ the {\em singular set\/}
 and denote by 
 $\Sigma_f:=\{p\in M^2\,|\,\text{$p$ is a singular point of $f$}\}$.
 A singular point $p\in \Sigma_f$ is called {\it non-degenerate\/} if
 the derivative $d\lambda$ of the signed area density function  does not
 vanish at $p$.
 This condition does not depend on choice of coordinate systems.
\end{definition}
It is well-known that a front can be considered locally
as a projection of a Legendrian immersion $L\colon{}U^2\to P(T^*N^3)$,
where $U^2$ is a domain in $\R^2$ and $P(T^*N^3)$ is the projective
cotangent bundle.
The canonical contact structure of the unit cotangent bundle 
$T^*_1N^3$ is the pull-back of that of $P(T^*N^3)$.
Since the contact structure on $P(T^*N^3)$ does not depend on the Riemannian
metric, the definition of front does not depend 
on the choice of the Riemannian metric $g$ and is invariant under 
diffeomorphisms of  $N^3$. 

\begin{definition}\label{def:limittangent}
 Let $f\colon{}M^2\to N^3$ be a front and $TN^3|_M$ the
 restriction of the tangent bundle of $N^3$ to $M^2$. 
 The subbundle $\E$ of rank $2$ on $M^2$ that
 is perpendicular 
 to the unit normal vector field $\nu$ of $f$ is called
 the {\it limiting tangent bundle\/} with respect to $f$.
\end{definition}
There exists a canonical vector bundle homomorphism
\[
  \psi \colon{}
    TM^2\ni X \longmapsto f_*X \in \E.
\]
The non-degenerateness in Definition~\ref{def:nondeg} 
is also independent of the choice of $g$
and can be described in terms of the limiting tangent bundle:

\begin{proposition}\label{prop:limit}
 Let $f\colon{}U\to N^3$ be a front defined on a domain $U$ in $\R^2$
 and $\E$ the limiting tangent bundle. 
 Let  $\mu\colon{}(U;u,v)\to \E^* \wedge \E^*$ 
 be an arbitrary fixed nowhere vanishing section.
 Then a singular point $p\in M^2$ is non-degenerate if and only
 if the derivative $dh$ of the function
 $h:=\mu\bigl(\psi(\partial/\partial u),\psi(\partial/\partial v)\bigr)$
  does not vanish at $p$.
\end{proposition}

\begin{proof}
Let $\mu_0$ be the $2$-form that is the restriction of
the $2$-form $\iota_\nu\mu_g$ to $M^2$, 
where $\iota_\nu$ denotes the interior product
and $\mu_g$ is the volume element of $g$.  
Then $\mu_0$ is a nowhere vanishing section on $\E^* \wedge \E^*$,
and the local signed area density function $\lambda$
is given by 
$\lambda=\mu_0(\psi(\partial/\partial u),\psi(\partial/\partial v))$.

On the other hand, 
let  $\mu\colon{}(U;u,v)\to \E^* \wedge \E^*$ be an arbitrary fixed
nowhere vanishing section.
Then there exists a smooth function 
$\tau\colon{}U\to \R\setminus\{0\}$ such that 
$\mu=\tau \cdot \mu_0$ (namely $h=\tau\lambda$) and
\[
  dh(p)=d\tau(p) \cdot \lambda(p)+\tau(p) 
             \cdot d\lambda(p)=\tau(p) \cdot d\lambda(p) ,
\]
since $\lambda(p)=0$ for each singular point $p$.
Then $dh$ vanishes if and only if 
$d\lambda$ does as well.
\end{proof}

\begin{remark}\label{rem:newadd}
 A $C^\infty$-map $f:U^2\to M^3$ is called a {\it frontal\/}
 if it is a projection of isotropic map
 $L:U^2\to T^*_1M^3$,
 that is, the pull-back of the canonical contact form of
 $T_1N^3$
 by $L$ vanishes on $M^2$.
 The definition of non-degenerate singular points  
 and the above lemma do not use the properties 
 that $L$ is an immersion.
 So they hold for any frontals.
\end{remark}

Let $p\in M^2$ be a non-degenerate singular point.
Then by the implicit function theorem, the singular set near $p$
consists of a regular curve in the domain of $M^2$.
This curve is called the {\em singular curve\/} at $p$.
We denote the singular curve by 
\[
   \gamma\colon{}(-\varepsilon,\varepsilon)
          \ni t \longmapsto \gamma(t)\in M^2
     \qquad (\gamma(0)=p).
\]
For each $t\in (-\varepsilon,\varepsilon)$, there exists
a $1$-dimensional linear subspace of $T_{\gamma(t)}M^2$, 
called the {\it null direction\/}, which is
the kernel of the differential map $f_*$.
A non-zero vector belonging to the null direction is called a 
{\em null vector}.
One can choose a smooth vector field $\eta(t)$ along $\gamma(t)$
such that $\eta(t)\in T_{\gamma(t)}M^2$ is a null vector for 
each $t$, which is called a {\em null vector field}.
The tangential $1$-dimensional vector space of the
singular curve $\gamma(t)$ is called the {\em singular direction}.
\begin{fact}[Criteria for cuspidal edges and swallowtails \cite{KRSUY}]
\label{fact:intrinsic-criterion}
 Let $p$ be a non-degenerate singular point of a front $f$,
 $\gamma$ the singular curve passing through $p$, and
 $\eta$ a null vector field along $\gamma$.
 Then
 \begin{enumerate}
  \item[{\rm (a)}]\label{item:intrinsic-criterion-cuspidal} 
        $p=\gamma(t_0)$ is a cuspidal edge 
        {\rm(}that is,  $f$ is locally diffeomorphic to 
        $f_C$ of \eqref{eq:cuspidal-swallow} in the introduction{\rm)}
	if and only if the null direction and the singular direction
	are transversal, that is, $\det\bigl(\gamma'(t),\eta(t)\bigr)$
	does not vanish at $t=t_0$, where 
	$\det$ denotes the determinant of $2\times 2$ matrices
	and where
	we identify the tangent space in $T_{\gamma(t_0)}M^2$ with
	$\R^2$.
  \item[{\rm (b)}]\label{item:intrinsic-criterion-swallow} 
        $p=\gamma(t_0)$ is a swallowtail
        {\rm (}that is, $f$ is locally diffeomorphic to $f_S$ of 
        \eqref{eq:cuspidal-swallow}
        in the introduction{\rm )} 
	if and only if 
	\[
	  \det\bigl(\gamma'(t_0),\eta(t_0)\bigr)=0 \qquad\text{and}\qquad
	  \left.\frac{d}{dt}\right|_{t=t_0}\!\!
	  \det\bigl(\gamma'(t),\eta(t)\bigr)\neq 0
	\]
	hold.
 \end{enumerate}
\end{fact}

For later computation, it is convenient to take a local coordinate
system $(u,v)$ centered at a given non-degenerate singular point 
$p\in M^2$ as follows:
\begin{itemize}
 \item the coordinate system $(u,v)$ is compatible with the orientation
        of $M^2$,
 \item the $u$-axis is the singular curve, and
 \item there are no singular points other than the $u$-axis.
\end{itemize}
\noindent
We call such a coordinate system $(u,v)$ an 
{\em adapted coordinate system\/} with respect to $p$.
In these coordinates, the signed area density function $\lambda(u,v)$ 
vanishes on the $u$-axis.
Since $d\lambda\ne 0$, $\lambda_v$ never vanishes on the $u$-axis.
This implies that 
\begin{equation}\label{eq:area-density-sign}
  \text{the signed area density function $\lambda$ changes sign on singular
   curves,}
\end{equation}
that is, the singular curve belongs to the boundary of $M_+$ and $M_-$.

Now we suppose that a singular curve $\gamma(t)$ on $M^2$
consists of cuspidal edges. 
Then we can choose the null vector fields $\eta(t)$ such that
$\bigl(\gamma'(t),\eta(t)\bigr)$ is a positively oriented frame field
along $\gamma$.
We then define 
the {\em singular curvature function\/} along $\gamma(t)$
as follows:
\begin{equation}\label{eq:def-singular-curvature}
 \kappa_s(t):=\sign\bigl(d\lambda(\eta)\bigr)\,
  \frac{\mu_g\bigl(\hat \gamma'(t),\hat \gamma''(t),\nu\bigr)}
  {|\hat \gamma'(t)|^3}.
\end{equation}
Here, we denote 
$|\hat\gamma'(t)|=g\bigl(\hat\gamma'(t),\hat\gamma'(t)\bigr)^{1/2}$,
\begin{equation}\label{eq:singular-curve-image}
    \hat \gamma(t)=f(\gamma(t)),\qquad
    \hat \gamma'(t)=\frac{d\hat \gamma(t)}{dt},\quad
    \text{and}\quad
    \hat \gamma''(t)=D_{t}\hat \gamma'(t),
\end{equation}
where $D$ is the Levi-Civita connection 
and $\mu_g$ the volume element of $(N^3,g)$.

\begin{figure}
 \begin{center}
  \input{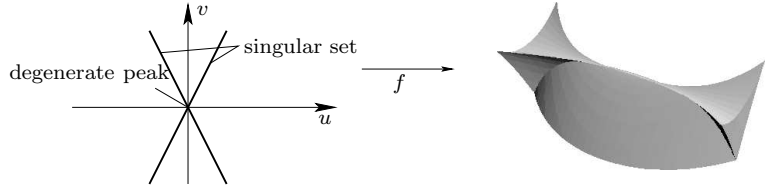}
 \end{center}
 \caption{A double swallowtail (Example~\ref{ex:degenerate-peak}).}
 \label{fig:double-swallow}
\end{figure}

We take an adapted coordinate system $(u,v)$ and write
the null vector field $\eta(t)$ as
\begin{equation}\label{eq:null-vector-adapted}
 \eta(t)=
  a(t)\frac{\partial}{\partial u}+
  e(t)\frac{\partial}{\partial v},
\end{equation}
where $a(t)$ and $e(t)$ are $C^\infty$-functions.
Since $(\gamma',\eta)$ is a positive frame, we have $e(t)>0$. 
Here, 
\begin{equation}\label{eq:non-degenerate-adapted}
 \lambda_u=0\qquad\text{and}\qquad  \lambda_v\neq 0
  \qquad \text{(on the $u$-axis)}
\end{equation}
hold,
and then  $d\lambda\bigl(\eta(t)\bigr)=e(t)\lambda_v$.
In particular, we have
\begin{equation}\label{eq:singular-curvature-sign-adapted}
 \sign\bigl(d\lambda(\eta)\bigr)=
  \sign(\lambda_v)
  =
  \begin{cases}
   +1 & \text{if the left-hand side of $\gamma$ is $M_+$}, \\
   -1 & \text{if the left-hand side of $\gamma$ is $M_-$}.
  \end{cases}
\end{equation}
So we have the following expression:
in an adapted coordinate system $(u,v)$,
\begin{equation}\label{eq:singular-curvature-adapted}
 \kappa_s(u):=\sign(\lambda_v)
     \frac{\mu_g(f_u,f_{uu},\nu)}
     {|f_u|^3},
\end{equation}
where $f_{uu}=D_uf_u$
and $|f_u|=g(f_u,f_u)^{1/2}$.
\begin{theorem}[Invariance of the singular curvature]
\label{thm:invariance-singular-curvature}
 The definition \eqref{eq:def-singular-curvature} of the singular
 curvature does not depend on the parameter $t$, nor the orientation of
 $M^2$,  nor the choice of $\nu$,
 nor the orientation of the singular curve.
\end{theorem}
\begin{proof}
 If the orientation of $M^2$ reverses, then $\lambda$ and $\eta$
 both change sign.
 If $\nu$ is changed to $-\nu$, so does $\lambda$.
 If $\gamma$ changes orientation, both $\gamma'$ and $\eta$ 
 change sign. 
 In all cases, 
 the sign of $\kappa_s$ is unchanged.
\end{proof}
\begin{remark}\label{rem:geodesic-curvature}
 We have the following expression
 \begin{multline*}
 \kappa_s=\sign\bigl(d\lambda(\eta)\bigr)\,
     \frac{
        \mu_0(\hat \gamma'',\nu,\hat \gamma'/|\hat \gamma'|)}{%
        |\hat \gamma'|^2}
    =\sign\bigl(d\lambda(\eta)\bigr)\,
      \frac{g(\hat \gamma'',n)}{|\hat \gamma'|^2} \quad
      \left(
         n:=\nu\times_g \frac{\hat\gamma'}{|\hat \gamma'|}
      \right).
  \end{multline*}
 Here,
 the vector product operation $\times_g$ in $T_xN^3$ is defined by 
 $a\times_g b:=*(a\wedge b)$, 
 under the identification 
 $TN^3\ni X \leftrightarrow  g(X,~)\in T^*N^3$,
 where $*$ is the Hodge $*$-operator.
 If $\gamma(t)$ is not a singular curve,
 $n(t)$ is just the conormal vector of $\gamma$.
 We call $n(t)$ the 
 {\em limiting conormal vector\/}, and
 $\kappa_s(t)$ can be considered as the limiting geodesic curvature
 of (regular) curves with the singular curve on their right-hand sides.
\end{remark}
\begin{proposition}[Intrinsic formula for the singular curvature]
\label{prop:intrinsic-singular-curvature}
 Let $p$ be a point of a cuspidal edge of a front $f$, and 
 $(u,v)$ an adapted coordinate system at $p$
 such that $\partial/\partial v$ gives the null direction. 
 Then the singular curvature is given by
 \[
   \kappa_s(u)=\frac{-F_vE_u+2EF_{uv}-E E_{vv}}{E^{3/2}\lambda_v},
 \]
 where $E=g(f_u,f_u),F=g(f_u,f_v),G=g(f_v,f_v)$,
 and where $\lambda$ is the signed area density function 
 with respect to $(u,v)$.
\end{proposition}
\begin{proof}
 Fix $v>0$ and denote by $\gamma(u)=(u,v)$ the $u$-curve.
 Then the unit vector
  \[
    n(u)=\frac{1}{\sqrt{E}\sqrt{EG-F^2}}
        \left(-F \frac{\partial}{\partial u}+
               E \frac{\partial}{\partial v}\right)
 \]
 gives the conormal vector such that
  $\bigl(\gamma'(u),n(u)\bigr)$ is a positive frame.
 Let $\nabla$ be the Levi-Civita connection on $\{v>0\}$
 with respect to the induced metric $ds^2=Edu^2+2Fdudv+Gdv^2$,
 and $s$ the arclength parameter of $\gamma(u)$.
 Then we have
 \[
      \nabla_{\gamma'(s)}^{}\gamma'(s)
       =
      \frac{1}{\sqrt{E}}
      \nabla_{\partial/\partial u}^{}
          \left( \frac{1}{\sqrt{E}}\frac{\partial}{\partial u}\right)
         \equiv \frac{\Gamma_{11}^2}{E} 
         \frac{\partial}{\partial v}
         \mod \frac{\partial}{\partial u},
 \]
 where $\Gamma_{11}^2$ is the Christoffel symbol given by
 \[
   \Gamma_{11}^2=\frac{-F E_{u}+2E F_u-EE_v}{2(EG-F^2)}.
 \]
 Since $\lambda^2=EG-F^2$  and $g(f_u,n)=0$,
 the geodesic curvature of $\gamma$ is given by
\[
    \kappa_g=g\bigl(\nabla_{\gamma'(s)},\gamma'(s),{n(s)}\bigr)
     =\frac{\sqrt{EG-F^2}\, \Gamma_{11}^2}{E^{3/2}}
     =\frac{-F E_{u}+2E F_u-EE_v}{|\lambda| E^{3/2}}.
\]
 Hence, by Remark~\ref{rem:geodesic-curvature}, the singular curve of
 the $u$-axis is
 \[
   \kappa_s=\sign(\lambda_v)\lim_{v\to 0}\kappa_g
    =\sign(\lambda_v)
     \lim_{v\to 0}
       \frac{-F E_{u}+2E F_u-EE_v}{|\lambda| E^{3/2}}.
 \]
 It is clear that all of
 $\lambda$, $F$ and $F_u$ tend to zero as $v\to 0$.
 Moreover, we have
 \[
     E_v
     =2g(D_vf_u,f_u)=2g(D_uf_v,f_u)
     =2\frac{\partial}{\partial v}g(f_v,f_u)-2g(f_v,D_uf_u)\to 0
 \]
 as $v\to 0$,  and the right differential $|\lambda|_v$ is equal to 
 $|\lambda_v|$ since $\lambda(u,0)=0$. 
 By L'Hospital's rule, we have
 \[
    \kappa_s
         =\sign(\lambda_v)
          \frac{-F_v E_{u}+2E F_{uv}-EE_v}{|\lambda|_v E^{3/2}}
         =\frac{-F_v E_{u}+2E F_{uv}-EE_v}{\lambda_v E^{3/2}},
 \]
 which is the desired conclusion.
\end{proof}
\begin{example}[Cuspidal parabolas]
\label{ex:parabola}
 Define a map $f$ from $\R^2$ to the Euclidean $3$-space $(\R^3,g_0)$ as
 \begin{equation}\label{eq:parabola}
    f(u,v)=(au^2+v^2,bv^2+v^3,u)
            \qquad (a,b\in \R).
 \end{equation}
 Then we have
 $f_u=(2 a u, 0, 1)$, $f_v=(2v, 2bv + 3v^2, 0)$.
 This implies that the $u$-axis is the singular curve, 
 and the $v$-direction is the null direction.
 The unit normal vector 
 and the signed area density $\lambda=\mu_{g_0}(f_u,f_v,\nu)$ are given
 by
 \begin{multline}\label{eq:parabola-normal}
    \nu = \frac{1}{\delta}
         \bigl(-3v-2b,2,2au(3v+2b)\bigr),\qquad
    \lambda = v\delta, \\
   \text{where}\qquad
   \delta = \sqrt{4+(1+4a^2u^2)(4b^2+12bv+9v^2)}.
 \end{multline}
 In particular, since $d\nu(\partial/\partial v)=\nu_v\neq 0$ on
 the $u$-axis, $(f,\nu)\colon{}\R^2\to \R^3\times S^2=T_1\R^3$
 is an immersion, i.e.\ $f$ is a front, and
 each point of the $u$-axis is a cuspidal edge. 
 The singular curvature is given by
 \begin{equation}\label{eq:curvature-parabola}
  \kappa_s(u)=
     \frac{2a}{(1+4a^2u^2)^{3/2}\sqrt{1+b^2(1+4a^2u^2)}}.
 \end{equation}
 When $a>0$ (resp.\ $a<0$), that is,
 the singular curvature is positive (resp.\ negative), 
 we shall call $f$ a {\em cuspidal elliptic\/}
 (resp.\ {\em hyperbolic}) {\em parabola\/}
 since the figure looks like a elliptic 
 (resp.\ hyperbolic) parabola,
 as seen in Figure~\ref{fig:parabola} in the introduction.
\end{example}
\begin{definition}[Peaks]\label{def:peak}
 A singular point $p\in M^2$ (which is not a cuspidal edge)
 is called a {\em peak\/} 
 if there exists a coordinate neighborhood $(U;u,v)$ of $p$
 such that
 \begin{enumerate}
  \item\label{item:peak-1} 
       there are no singular points other than cuspidal edges on 
       $U\setminus \{p\}$,
  \item\label{item:peak-2} 
       the rank of the derivative $f_*\colon{}T_pM^2\to T_{f(p)}N^3$
       at $p$ is equal to $1$, and
  \item\label{item:peak-3} 
       The singular set of $U$ consists of finitely many regular
       $C^1$-curves starting at $p$.
       The number $2m(p)$ of these curves
       is called the {\em number of cuspidal edges\/} starting at $p$.
 \end{enumerate}
 If a peak is a non-degenerate singular point, it
 is called a {\em non-degenerate peak}.
\end{definition}
Swallowtails are  examples of non-degenerate peaks.
A front which admits cuspidal edges and peaks is called
{\em a front which admits at most peaks}.
There are degenerate singular points which are not peaks.
Typical examples are cone-like singularities 
which appear in rotationally symmetric surfaces in $\R^3$ of positive
constant Gaussian curvature. 
However, since  generic fronts (in the local sense) have only cuspidal edges
and swallowtails,
the set of fronts which admits at most peaks covers a sufficiently wide
class of fronts.
\begin{example}[A double swallowtail]
\label{ex:degenerate-peak}
 Define a map $f\colon{}\R^2\to\R^3$ as 
 \[
     f(u,v) := (2u^3-uv^2,3u^4-u^2v^2,v).
 \]
 Then 
 \[
    \nu = \frac{1}{\sqrt{1+4u^2(1+u^2v^2)}}(-2u,1,-2u^2v)
 \]
 is the unit normal vector to $f$.
 The pull-back of the canonical metric of $T_1\R^3=\R^3\times S^2$
 by $(f,\nu)\colon{}\R^2\to\R^3\times S^2$ is positive definite.
 Hence $f$ is a front.
 The signed area density function is 
 $\lambda = (v^2-6u^2)\sqrt{1+4u^2(1+u^2v^2)}$, and 
 then the singular set is 
 $\Sigma_f=\{v=\sqrt{6}u\}\cup \{v=-\sqrt{6}u\}$.
 In particular, $d\lambda =0$ at $(0,0)$.
 The first fundamental form of $f$ is expressed as $ds^2=dv^2$ at
 the origin, 
 which is of rank one.
 Hence the origin is a degenerate peak
 (see Figure~\ref{fig:double-swallow}).
\end{example}

To analyze the behavior of the singular curvature near a peak,
we prepare the following proposition.
\begin{proposition}%
[Boundedness of the singular curvature measure]
\label{prop:bounded-curvature-measure}
 Let $f\colon{}M^2\to (N^3,g)$ be a front with a peak  $p$.
 Take $\gamma\colon{}[0,\varepsilon)\to M^2$ a singular curve of $f$
 starting from the singular point $p$. 
 Then $\gamma(t)$ is a cuspidal edge for $t>0$, and
 the singular curvature measure $\kappa_s\,ds$ is continuous
 on $[0,\varepsilon)$, where $ds$ is the arclength-measure.
 In particular, the limiting tangent vector
 $\displaystyle\lim_{t\to 0}\hat\gamma'(t)/|\hat\gamma'(t)|$ 
 exists, 
 where $\hat\gamma=f\circ \gamma$.
\end{proposition}

\begin{proof}
 Let $ds^2$ be the first fundamental form of $f$.
 Since $p$ is a peak, the rank $ds^2$ is $1$ at $p$ 
 and then one of the eigenvalues is $0$ and the other is not.
 Hence the eigenvalues of $ds^2$ are of multiplicity one 
 on a neighborhood of $p$.
 Hence one can choose a local coordinate system $(u,v)$ around $p$
 such that each coordinate curve is tangent to an eigendirection of 
 $ds^2$.
 In particular, we can choose $(u,v)$ such that $\partial/\partial v$ 
 is the null vector field on $\gamma$.
  In such a coordinate system,  $f_v=0$ and $D_t f_v=0$ hold on
 $\gamma$.
 Then the derivatives of $\hat\gamma=f\circ\gamma$ are
 \[
    \hat\gamma' = u'f_u,\qquad
    D_t\hat\gamma' = u'' f_u + u'D_t f_u\qquad
    \left('=\frac{d}{dt}\right),
 \]
 where $\gamma(t)=\bigl(u(t),v(t)\bigr)$.
 Hence 
 \begin{equation}\label{eq:singular-curvature-peak}
    \kappa_s =
    \pm \frac{\mu_g(\hat\gamma',D_t\hat\gamma',\nu)}{|\hat\gamma'|^3}
     = \pm \frac{\mu_g(f_u,D_t f_u,\nu)}{|u'|\,|f_u|^3},
 \end{equation}
 where $|X|^2=g(X,X)$ for $X\in TN^3$.
 Since $ds=|\hat\gamma'|\,dt=|u'|\,|f_u|\,dt$ and $f_u\neq 0$,
 \[
    \kappa_s\,ds =\pm \frac{\mu_g(f_u,D_tf_u,\nu)}{|f_u|^2}\,dt
 \]
 is bounded.
\end{proof}

\begin{figure}\footnotesize
 \begin{center}
  \begin{tabular}{c@{\hspace{2cm}}c}
   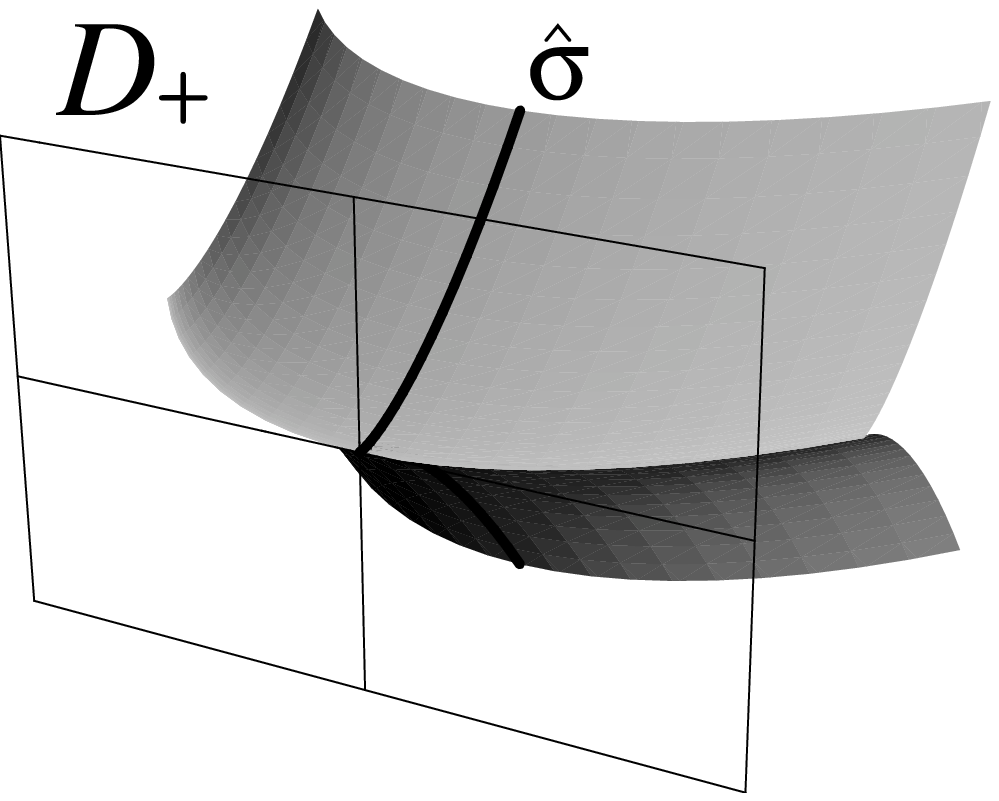&
   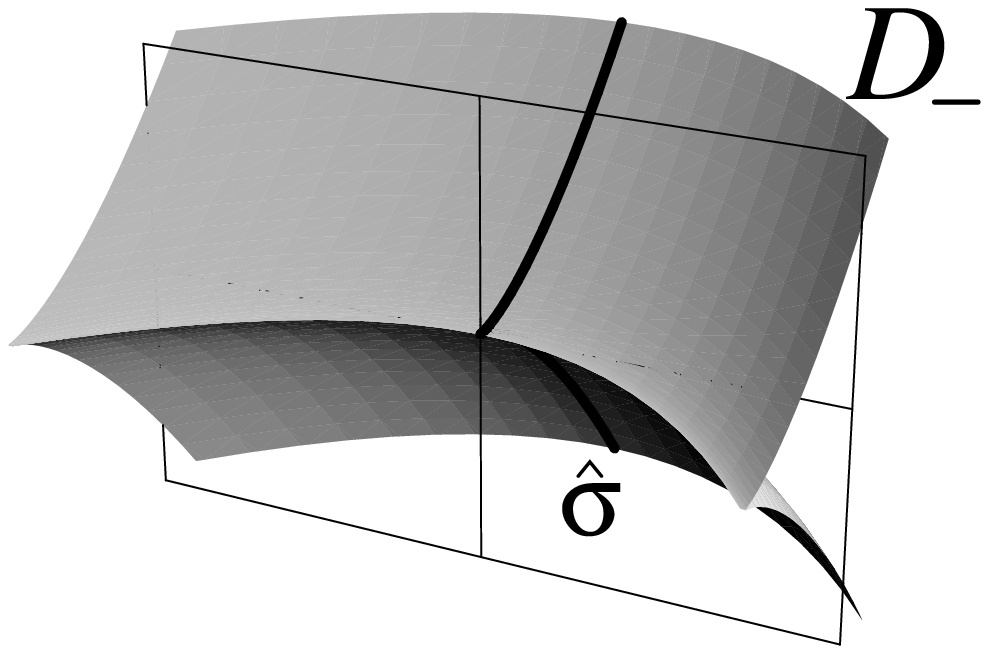\\
   positively curved & negatively curved
  \end{tabular}\\
  The principal half-spaces are behind the rectifying plane.
 \end{center}
 \caption{The principal half-spaces of cuspidal edges.}
 \label{fig:principal-space}
\end{figure}

To analyze the behavior of the singular curvature near 
a non-degenerate peak, 
we give another expression of the singular curvature measure:
\begin{proposition}\label{prop:bdd2}
 Let $(u,v)$ be an adapted coordinate system of $M^2$.
 Suppose that $(u,v)=(0,0)$ is a non-degenerate peak.
 Then the singular curvature measure has the expression
 \begin{equation}\label{eq:sing-3}
  \kappa_s(u)ds=\sign(\lambda_v)\frac{\mu_g(f_v,f_{uv},\nu)}{|f_v|^2}du,
 \end{equation}
 where $ds$ is the arclength-measure and
 $f_{uv}:=D_uf_v=D_vf_u$.
 In particular, the singular curvature measure
 is smooth along the singular curve.
\end{proposition}
\begin{proof}
 We can take the null direction
 $\eta(u)=a(u)(\partial/\partial u)+ e(u)(\partial/\partial v)$
 as in \eqref{eq:null-vector-adapted}.
 Since the peak is not a cuspidal edge, 
 $\eta(0)$ must be proportional to $\partial_u$.
 In particular, we can multiply 
 $\eta(u)$ by a non-vanishing function
 and may assume that $a(u)=1$.
 Then $f_u+e(u)f_v=0$ and by differentiation
 we have
 $f_{uu}+e_uf_v+e f_{uv}=0$,
 that is,
 \[
   f_u=-e f_v,\quad 
   f_{uu}=-e_u f_v-e f_{uv}.
 \]
 Substituting them into \eqref{eq:singular-curvature-adapted}, we have
 \eqref{eq:sing-3} using the relation
 $ds=|\hat \gamma'|dt=|f_u|dt$.
\end{proof}
\begin{corollary}%
[Behavior of the singular curvature near a non-degenerate peak]
\label{cor:singular-curvature-peak}\ 
 At a non-degenerate peak, 
 the singular curvature diverges to  $-\infty$.
\end{corollary}
\begin{proof}
 We take an adapted coordinate $(u,v)$ centered at the peak.
 Then
 \[
   \kappa_s(u) =\sign(\lambda_v)
              \frac{\mu_g(f_{v},f_{uv},\nu)}{|e(u)|\, |f_v|^3}.
 \]
 On the other hand,
 \[
   \mu_g(f_{v},f_{uv},\nu)=
      \mu_g(f_{v},f_{u},\nu)_v-\mu_g(f_{vv},f_{u},\nu)
   =
   (-\lambda)_v-\mu_g(f_{vv},f_{u},\nu).
 \]
 Since $f_u(0,0)=0$ we have
 \[
   \left.\sign(\lambda_v)
   \frac{\mu_g(f_{v},f_{uv},\nu)}{|f_v|^3}\right|_{(u,v)=(0,0)}
   \!\!=
   -\frac{|\lambda_v(0,0)|\hphantom{^3}}{|f_v(0,0)|^3}<0.
 \]
 Since $e(u)\to 0$ as $u\to 0$, we have the assertion.
\end{proof}
\begin{example}[The discriminant set of $s^3+zs^2+ys+x$]
\label{ex:swallowtail}
 The typical example of peaks is a swallowtail.
 We shall compute the singular curvature of the swallowtail 
 $f(u,v)=(3u^4+u^2v ,4u^3+2uv,v)$ at $(u,v)=(0,0)$ given in 
 the introduction, 
 which is the discriminant set
 $\{(x,y,z)\in \R^3\,;\, 
  F(x,y,z,s)=F_{s}(x,y,z,s)=0~ 
  \text{for $s\in \R$}\}$
 of the polynomial $F:=s^3+zs^2+ys+x$ in $s$.
 Since
 $f_u\times f_v=2(6u^2+v)(1,-u,u^2)$, the singular curve is 
 $\gamma(t)=(t,-6t^2)$
 and the unit normal vector is given by
 $\nu=(1,-u,u^2)/\sqrt{1+u^2+u^4}$.
 We have 
 \[
   \kappa_s(t)=
        \frac{\det(\hat\gamma',\hat\gamma'',\nu)}{|\hat\gamma'|^3}
              =-\frac{\sqrt{1+t^2+t^4}}{6|t|(1+4t^2+t^4)^{3/2}},
 \]
 which shows the singular curvature tends to $-\infty$ when $t\to 0$.
\end{example}
\begin{definition}[Null curves]
\label{def:null}
 Let $f\colon{}M^2\to N^3$ be a front.
 A regular curve $\sigma(t)$ in $M^2$ is called a 
 {\em null curve\/} of $f$ if
 $\sigma'(t)$ is a null vector at each singular point.
 In fact,  $\hat \sigma(t)=f\bigl(\sigma(t)\bigr)$ 
 looks like the curve (virtually) transversal
 to the cuspidal edge,  in spite of
 $\hat \sigma'=0$, and $D_t\hat \sigma'$ gives the ``tangential''
 direction of the surface at the singular point.
\end{definition}
\begin{theorem}[A geometric meaning for the singular curvature]
\label{thm:PN-cuspidal-edge}
 Let $p$ be a cuspidal edge,
 $\gamma(t)$ a singular curve 
 parametrized by the arclength $t$
 with $\gamma(0)=p$,
 and $\sigma(s)$  a null curve passing through  $p=\sigma(0)$.
 Then the sign of
 \[
    g\bigl(\ddot{\hat \sigma}(0),\hat \gamma''(0)\bigr)
 \]
 coincides with that of the singular curvature at $p$,
 where $\hat\sigma=f(\sigma)$, $\hat\gamma=f(\gamma)$,
 \[
   \dot{\hat\sigma}=\frac{d\hat\sigma}{ds},\qquad
    \hat\gamma'=\frac{d\hat\gamma}{dt},\qquad
    \ddot{\hat\sigma}=D_s\left(\frac{d\hat\sigma}{ds}\right),
   \quad\text{and}\quad
    \hat\gamma''=D_t\left(\frac{d\hat\gamma}{dt}\right).
 \]
\end{theorem}
\begin{proof}
 We can take an adapted coordinate system
 $(u,v)$ around $p$ such that $\eta:=\partial/\partial v$
 is a null vector field on the $u$-axis.
 Then
 $f_v=f_*\eta$ vanishes on the $u$-axis, and it holds that
 $f_{uv}:=D_v f_u=D_u f_v=0$ on the $u$-axis.
 Since the $u$-axis is parametrized by the arclength, we have
 \begin{equation}
  \label{eq:u-axis-unit-length}
  g(f _{uu},f_u)=0 \qquad \text{on the $u$-axis}
   \qquad\left(f_{uu}=D_u f_u\right).
 \end{equation}
 Now let $\sigma(s)=\bigl(u(s),v(s)\bigr)$
 be a null curve such that $\sigma(0)=(0,0)$.
 Since $\dot\sigma(0)$ is a null vector, $\dot u(0)=0$, where
 $\dot{~}=d/ds$.
 Moreover, since $f_v(0,0)=0$ and $f_{uv}(0,0)=0$, we have  
 \begin{align*}
  \ddot{\hat \sigma}(0)
    &= D_s (\dot u f_u+ \dot v f_v) =  
    \ddot u f_u+  \ddot v f_v+  \dot u^2 D_{u} f_u +
            2\dot  u \dot v D_uf_v +\dot v^2D_v f_v \\
    &=  \ddot u f_u+  \dot v^2 D_v f_v
     =  \ddot u f_u(0,0)+ \dot v^2 f_{vv}(0,0),
 \end{align*}
 and by \eqref{eq:u-axis-unit-length},
 \[
   g\bigl(\ddot{\hat\sigma}(0)),\hat\gamma''(0)\bigr)
    =g\bigl(f_{uu}(0,0),\ddot u f_u+ \dot v^2 f_{vv}(0,0)\bigr)
    =\dot v^2 g\bigl(f_{uu}(0,0),f_{vv}(0,0)\bigr).
 \]
 Now we can write
 $f_{vv}=a f_u+b(f_u\times_g \nu)+c \nu$, where $a,b,c\in \R$.
 Then
 \begin{align*}
   c&=g(f_{vv},\nu)=g(f_v,\nu)_v-g(f_v,\nu_{v})=0,\\
   b&=g(f_{vv}, f_u\times_g \nu)=
     g(f_v, f_u\times_g \nu)_v=-\lambda_v,
 \end{align*}
 where we apply the scalar triple product formula
 $g(X, Y\times_g Z)=\mu_g(X,Y,Z)$ for $X,Y,Z\in T_{f(0,0)}N^3$.
 Thus
 \begin{align*}
  g\bigl(\ddot{\hat\sigma}(0),\hat\gamma''(0)\bigr)
    &=\dot v^2 g\bigl(f_{uu}, a f_u-\lambda_v(f_u\times_g \nu)\bigr)
    =-\dot v^2 \lambda_v\,g(f_{uu}, f_u\times_g \nu)\\
    &=\dot v^2  \lambda_v\,\mu_g(\hat \gamma',\hat \gamma'',\nu)
    = \dot v^2|\lambda_v|\,\kappa_s(0).
 \end{align*}
 This proves the assertion.
\end{proof}
In the case of fronts in the Euclidean $3$-space $\R^3=(\R^3,g_0)$,
positively curved cuspidal edges and 
negatively curved cuspidal edges 
look like cuspidal elliptic parabola or 
hyperbolic parabola 
(see Example~\ref{ex:parabola} and Figure~\ref{fig:parabola}),
respectively.
More precisely, we have the following:
\begin{corollary}%
\label{cor:positive-negative-cuspidal-edge}
 Let $f\colon{}M^2\to(\R^3,g_0)$ be a front, $p\in M^2$ a 
 cuspidal edge point and $\gamma$  a singular curve with $\gamma(0)=p$.
 Let $T$ be the rectifying plane of the singular curve 
 $\hat\gamma=f\circ\gamma$ at $p$, 
 that is, the plane perpendicular to the principal normal vector of
 $\hat\gamma$.
 When the singular curvature at $p$ is positive
 {\rm (}resp.\ negative{\rm)}, every null
 curve $\sigma(s)$ passing through $\sigma(0)=p$ lies on the
 same side $D_+$ {\rm(}resp.\ the opposite side $D_-${\rm)}
 of the principal normal vector of $\hat\gamma$ at $p$
 for sufficiently small $s$.
 Moreover, if the singular curvature is positive, 
 the image of the neighborhood of $p$ itself lies in $D_+$
 {\rm (}see Figures \ref{fig:parabola}
    and \ref{fig:principal-space}{\rm )}.
\end{corollary}
\begin{definition}\label{def:principal}
 The half-space 
 in  Corollary~\ref{cor:positive-negative-cuspidal-edge}
 bounded by the rectifying plane of the singular curve 
 and in which
 the null curves lie is called the {\it principal half-space\/}
 at the cuspidal edge.
 The surface lies mostly in this half-space.
 When the singular curvature is positive, the surface is locally
 inside the principal half-space.
\end{definition}

\begin{proof}[Proof of Corollary~\ref{cor:positive-negative-cuspidal-edge}]
 Let $(u,v)$ be the same coordinate system at $p$ as in
 the proof of Proposition~\ref{thm:PN-cuspidal-edge} and assume 
 $f(0,0)=0$.
 Since $N^3=\R^3$, with $f_{uu}=\partial^2 f/\partial u^2$ etc.,
 we have the following Taylor expansion:
 \begin{equation}\label{eq:taylor}
   f(u,v) = f_u(0,0) + 
     \frac{1}{2}\bigl(f_{uu}(0,0)u^2 + f_{vv}(0,0)v^2\bigr)+
     o(u^2+v^2).
 \end{equation}
 Here, $u$ is the arclength parameter of $\hat\gamma(u)=f(u,0)$.
 Then $g_0(f_u,f_{uu})=0$ holds on the $u$-axis.  
 Thus 
 \[
    g_0\bigl(f(u,v),f_{uu}(0,0)\bigr)
       =\frac{1}{2}u^2|f_{uu}(0,0)|^2 +
        \frac{1}{2}v^2 g_0\bigl(f_{vv}(0,0),f_{uu}(0,0)\bigr)
           + o(u^2+v^2).
 \]
 If the singular curvature is positive,
 Theorem~\ref{thm:PN-cuspidal-edge} implies 
 $g_0\bigl(f(u,v),f_{uu}(0,0)\bigr)>0$ on a neighborhood of $p$.
 Since $f_{uu}(0,0)$ is the principal curvature vector of $\hat\gamma$
 at $p$, $f(u,v)$ lies in the same side of $T$ as the principal normal.
 
 Next we suppose that the singular curvature is negative at $p$.
 We can choose a coordinate system in which the null curve
 is written as $\sigma(v)=(0,v)$.
 Then by \eqref{eq:taylor} and Theorem~\ref{thm:PN-cuspidal-edge},
 \[
    g_0\bigl(f(0,v),f_{uu}(0,0)\bigr)
    = v^2 g_0\bigl(f_{vv}(0,0),f_{uu}(0,0)\bigr)+o(v^2)<0
 \]
 for sufficiently small $v$.
 Hence we have the conclusion.
\end{proof}
\begin{example}[Fronts with Chebyshev net]
\label{ex:chebyshev}
 A front $f\colon{}M^2\to \R^3$ is said to be of 
 {\em constant Gaussian curvature $-1$\/} 
 if the set $W=M^2\setminus\Sigma_f$ of regular points are dense in $M^2$
 and $f$ has constant Gaussian curvature $-1$ on $W$.
 Then $f$ is a projection of the Legendrian immersion
 $L_f\colon{}M^2\to T_1\R^3$, 
 and the pull-back $d\sigma^2=|df|^2+|d\nu|^2$ of the 
 Sasakian metric on $T_1\R^3$ by $L_f$ is flat. 
 Thus for each $p\in M^2$, there exists a coordinate neighborhood
 $(U;u,v)$ such that $d\sigma^2=2(du^2+dv^2)$.
 The two  different families of asymptotic curves on $W$ are 
 all geodesics of $d\sigma^2$, giving two foliations of $W$.
 Moreover, they are mutually orthogonal with respect to $d\sigma^2$.
 Then one can choose the $u$-curves and $v$-curves to all be
 asymptotic curves on $W\cap U$.
 For such a coordinate system $(u,v)$,
 the first and second fundamental forms are
 \begin{equation}\label{eq:chebyshev-front}
   ds^2 = du^2 + 2\cos\theta\,du\,dv + dv^2 ,\qquad
   h= 2\sin\theta\,du\,dv,
 \end{equation}
 where $\theta=\theta(u,v)$ is the angle between the two asymptotic curves.
 The coordinate system $(u,v)$ as in \eqref{eq:chebyshev-front}
 is called the {\em asymptotic Chebyshev net\/} around $p$.
 The sine-Gordon equation $\theta_{uv}=\sin\theta$ is the
 integrability condition of \eqref{eq:chebyshev-front}, that is, 
 if $\theta$ satisfies the sine-Gordon equation, then there exists a
 corresponding front $f=f(u,v)$.

 For such a front, we can choose the unit normal vector $\nu$
 such that $f_u \times f_v = \sin\theta \,\nu$ holds,
 that is, $\lambda=\sin\theta$.
 The singular sets are characterized by $\theta\in \pi\Z$.
 We write $\varepsilon = e^{\pi i \theta}=\pm 1$ at a singular point.
 A given singular point is non-degenerate if and only if  $d\theta\ne 0$.
 Moreover, the cuspidal edges are characterized by 
 $\theta_u-\varepsilon\theta_v\ne 0$, and
 the swallowtails are characterized by
 $\theta_u+\varepsilon\theta_v\ne 0$, 
 $\theta_u-\varepsilon\theta_v= 0$ and
 $\theta_{uu}+\theta_{vv}\ne 0$.
 By a straightforward calculation applying 
 Proposition~\ref{prop:intrinsic-singular-curvature},
 we have
 \[
   \kappa_s=-\varepsilon\frac{\theta_u\theta_v}{%
     |\theta_u-\varepsilon \theta_v|}
      \qquad (\varepsilon=e^{\pi i \theta}).
 \]
 Recently Ishikawa-Machida \cite{IM} showed that
 the generic singularities of such fronts
 are cuspidal edges or
 swallowtails,
 as an application of Fact~\ref{fact:intrinsic-criterion}. 
\end{example}
\section{The Gauss-Bonnet theorem}
\label{sec:GB}
In this section, we shall generalize the two types of 
Gauss-Bonnet formulas mentioned in the introduction to
compact fronts which admit at most peaks.
\begin{proposition}
 \label{thm:gaussian-curvature-extended}
 Let $f\colon{}M^2\to (N^3,g)$ be a front, and $K$ the
 Gaussian curvature of
 $f$ which is defined on the set of regular points of $f$.
 Then  
 $K\,d\hat A$ can be  continuously extended as a globally defined
 $2$-form on $M^2$, where $d\hat A$ is the signed area form as in
 \eqref{eq:signed-area-form}.
\end{proposition}
\begin{proof}
 Let $(u,v)$ be a local coordinate system compatible to the 
 orientation of $M^2$, and $S=(S^i_j)$ the (matrix representation of)
 the shape operator of $f$ which is defined on the set 
 of regular points $M^2\setminus\Sigma_f$.
 That is, the Weingarten equation holds:
 \[
    \nu_u = - S^1_1 f_u - S^2_1 f_v,\quad
    \nu_v = - S^1_2 f_u - S^2_2 f_v,\quad
   \text{where }
    \nu_u=D_{u}\nu,~
    \nu_v=D_{v}\nu.
 \]
 Since the extrinsic curvature is defined as $K_{\ext}=\det S$,
 we have
 \[
    \mu_g(\nu_u,\nu_v,\nu)=
    (\det S)\,\mu_g(f_u,f_v,\nu)=
    K_{\ext}\,\lambda,
 \]
 where $\lambda$ is the signed area density.
 Thus,
 \[
    K_{\ext}\,d\hat A =
    K_{\ext}\,\lambda\,du\wedge dv = 
    \mu_g(\nu_u,\nu_v,\nu)\,du\wedge dv
 \]
 is a well-defined smooth $2$-form on $M^2$.
 
 By the Gauss equation, the Gaussian curvature
 $K$ satisfies
 \begin{equation}\label{eq:k-int-ext}
   K = c_{N^3}+K_{\ext}, 
 \end{equation}
 where $c_{N^3}$ is the sectional curvature  of $(N^3,g)$ with respect
 to the tangent plane.
 Since $f_*T_pM^2\subset T_{f(p)}N^3$ is 
 the orthogonal complement of the normal vector $\nu(p)^{\perp}$,
 the tangent plane is well-defined on all of  $M^2$.
 Thus $c_{N^3}$ is a smooth function, and 
 \[
    K\,d\hat A = c_{N^3}\,d\hat A + K_{\ext}\,d\hat A
 \]
 is a smooth $2$-form defined on $M^2$.
\end{proof}

\begin{remark}\label{rmk:abs}
 On the other hand, 
 \[
    K\,dA =
         \begin{cases}
	     \hphantom{-}K\,d\hat A\qquad & (\text{on $M_+$}),\\
	     -K\,d\hat A\qquad & (\text{on $M_-$})
	 \end{cases}
 \]
 is bounded, and extends continuously to the closure of $M_+$ 
 and also to the closure of $M_-$.
 (However, $K\,dA$ cannot be extended continuously 
 to all of $M^2$.)
\end{remark}
Now we suppose that $M^2$ is compact and 
$f\colon{}M^2\to \R^3$ is a front which admits at most 
peak singularities.
Then the singular set coincides with 
$\partial M_+=\partial M_-$,
and $\partial M_+$ and $\partial M_-$ are
piecewise $C^1$-differentiable because all 
singularities are at most peaks,
and the limiting tangent vector of each singular curve starting
at a peak exists by Proposition \ref{prop:bounded-curvature-measure}.

\begin{figure}
\footnotesize
\begin{center}
 \input{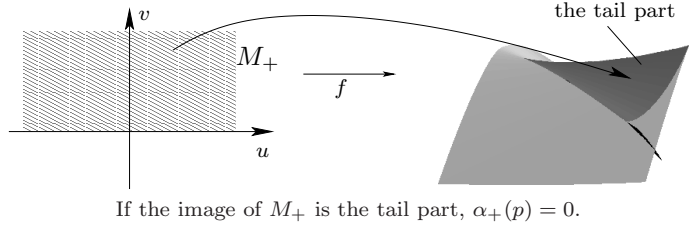}\\
 If the image of $M_+$ is the tail part, $\alpha_+(p)=0$.
\end{center}
\caption{A negative swallowtail.}
\label{fig:positive-negative-swallowtail}
\end{figure}

For a given peak $p$, 
let $\alpha_+(p)$ (resp.\ $\alpha_-(p)$)
be the sum of all the interior angles of $f(M_+)$ (resp.\ $f(M_-)$) at
$p$.
Then by definition, we have 
\begin{equation}
 \label{eq:sum-interior-angle}
  \alpha_+(p)+\alpha_-(p)=2\pi.
\end{equation}
Moreover, since the rank of $f_*$ is one at $p$, we have
(see \cite{SUY2})
\begin{equation}
 \label{eq:angle-range}
  \alpha_+(p),\,\,\alpha_-(p) \in \{0,\pi,2\pi\}.
\end{equation}
For example, $\alpha_+(p)=\alpha_-(p)=\pi$ when $p$ is 
a cuspidal edge.
If $p$ is a swallowtail,
$\alpha_+(p)=2\pi$ or $\alpha_-(p)=2\pi$.
If $\alpha_+(p)=2\pi$, $p$ is called a 
{\em positive swallowtail\/},
and is called a {\em negative swallowtail\/}
if $\alpha_-(p)=2\pi$
(see  Figure~\ref{fig:positive-negative-swallowtail}).
Since $K\,dA$, $K\,d\hat A$ and $\kappa_s\, ds$ are all bounded,
we get two Gauss-Bonnet formulas as follows:
\begin{theorem}[Gauss-Bonnet formulas for compact fronts]
 \label{thm:Gauss-Bonnet-compact}
 Let $M^2$ be a compact oriented $2$-manifold and
 $f\colon{}M^2\to (N^3,g)$ a front which admits at most peak
 singularities, and $\Sigma_f$ the singular set of $f$.
 Then 
 \begin{align}
  \int_{M^2}K\, dA+2\int_{\Sigma_f} \kappa_s\, ds&=
  2\pi \chi(M^2),
  \label{eq:Gauss-Bonnet-unsigned}\\
  \int_{M^2}K\, d\hat A-\sum_{p:\text{peak}} 
     \bigl(\alpha_+(p)-\alpha_-(p)\bigr)
  &=2\pi\bigl(\chi(M_+)-\chi(M_-) \bigr)
  \label{eq:Gauss-Bonnet-signed}
 \end{align}
 hold, where $ds$ is the arclength measure on the singular set.
\end{theorem}

\begin{remark}\label{rem:euler-characteristic}
 The integral $\int_{M^2}K\, d\hat A$
 is $2\pi$ times the Euler number 
 $\chi_{\E}^{}$  of the limiting tangent bundle $\E$
 (see \eqref{eq:Euler} in Section~\ref{sec:intrinsic}).
 When $N^3=\R^3$, $\chi_\E^{}/2$ is equal to the degree of the
 Gauss map.
\end{remark}
\begin{remark}\label{rem:generalization}
 These formulas are  generalizations of the two Gauss-Bonnet
 formulas in the introduction.
 If the surface is regular, the limiting tangent bundle 
 $\E$ coincides with the tangent bundle, 
 and the two Gauss-Bonnet formulas are the same.
\end{remark}%
\begin{proof}[Proof of Theorem~\ref{thm:Gauss-Bonnet-compact}]
 Although $\partial M_+$ and $\partial M_-$ are the same set, 
 their orientations are opposite.
 The singular curvature $\kappa_s$ 
 does not depend on the orientation of the singular curve
 and coincides with the limit of the geodesic curvature if we take the
 conormal vector in the positive direction with respect to the velocity
 vector of the singular curve.
 Thus we have 
 \begin{equation}\label{eq:integral-of-singular-curvature}
    \int_{\partial M_+}\!\! \kappa_s\, ds 
      +\int_{\partial M_-}\!\! \kappa_s\, ds=
    2\int_{\Sigma_f} \kappa_s\, ds.
 \end{equation}
 Then by the classical Gauss-Bonnet theorem, we have
 \begin{align*}
   2\pi\chi(M_+)
   &=\int_{M_+} K\,dA+\int_{\partial M_+}\!\!\kappa_s ds +
         \sum_{p:\text{peak}}\bigl(\pi m(p)-\alpha_+(p)\bigr),\\
   2\pi\chi(M_-)
   &=\int_{M_-} K\,dA+\int_{\partial M_-}\!\!\kappa_s\, ds +
         \sum_{p:\text{peak}}\bigl(\pi m(p)-\alpha_-(p)\bigr),
 \end{align*}
 where $2m(p)$ is the number of cuspidal edges starting at $p$
 (see Definition~\ref{def:peak}).
 Hence by \eqref{eq:integral-of-singular-curvature}, 
 \begin{align*}
  2\pi\chi(M^2)
&=  \int_{M^2} K\,dA +
                   2\int_{\Sigma_f}\kappa_s\,ds,\\
  2\pi\bigl(\chi(M_+)-\chi(M_-)\bigr) 
                   &= \int_{M^2} K\,d\hat A
                         -
                      \sum_{p:\text{peak}}
                      \bigl(\alpha_+(p)-\alpha_-(p)\bigr),
 \end{align*}
where we used 
\eqref{eq:sum-interior-angle} 
and 
$\chi(M^2)=\chi(M_+)+\chi(M_-)-
            \sum_{p:\text{peak}}\bigl(m(p)-1\bigr)$.
\end{proof}
We shall now define the completeness of fronts
and give Gauss-Bonnet formulas for non-compact fronts:
As defined in \cite{KUY2}, 
a front $f\colon{}M^2\to N^3$  is called {\it complete\/} 
if the singular set is compact and there exists a symmetric tensor $T$
with compact support such that $ds^2+T$ gives a complete Riemannian metric 
on $M^2$, where $ds^2$ is the first fundamental form of $f$.
On the other hand, as defined in \cite{KRSUY}, 
a front $f\colon{}M^2\to N^3$ is called {\em weakly complete\/} 
if the pull-back of the Sasakian metric of $T_1N^3$ by the Legendrian
lift $L_f\colon{}M^2\to T_1N^3$ is complete. 
Completeness implies weak completeness. 

Let $f\colon{}M^2\to N^3$ be a complete front with
finite absolute total curvature. 
Then there exists a compact $2$-manifold $\overline M^2$
without boundary and finitely many points  $p_1,\dots,p_k$ 
such that $M^2$ is diffeomorphic to 
$\overline M^2\setminus \{p_1,\dots,p_k\}$. 
We call the $p_i$'s the {\em ends\/} of the front $f$.
According to Theorem A of Shiohama \cite{S}, we define the limiting area
growth order
\begin{equation}\label{eq:area-growth}
  a(p_i)=\lim_{r\to \infty}
    \frac{\Area\bigl(B_0(r)\cap E_i\bigr)}{%
          \Area\bigl(B_{\R^2}(r)\cap E_i\bigr)},
\end{equation}
where $E_i$ is the punctured neighborhood of $p_i$ in $\overline{M}^2$.
\begin{theorem}[Gauss-Bonnet formulas for complete fronts]
\label{thm:Gauss-Bonnet-complete}
 Let $f\colon{}M^2\to (N^3,g)$ be a complete front with
 finite absolute total curvature, which has at most peak singularities,
 and write $M^2=\overline M^2\setminus\{p_1,\dots,p_k\}$.
 Then
 \begin{align}
  &\int_{M^2}\!\!K\, dA+2\int_{\Sigma_f} \kappa_s\, ds
    +\sum_{i=1}^ka(p_i)=2\pi\chi(M^2),
  \label{eq:Gauss-Bonnet-complete-unsigned}\\
  &\int_{M^2}\!\!K\, d\hat A-\sum_{p:\text{peak}} 
     \bigl(\alpha_+(p)-\alpha_-(p)\bigr)
       +\sum_{i=1}^k\varepsilon(p_i)a(p_i)
  =2\pi\bigl(\chi(M_+)-\chi(M_-)\bigr)
  \label{eq:Gauss-Bonnet-complete-signed}
 \end{align}
 hold, where $\varepsilon(p_i)=1$ 
 {\rm (}resp.\ $\varepsilon(p_i)=-1${\rm )}
 if the neighborhood $E_i$ of $p_i$ is contained in
 $M_{+}$ $($resp. $M_{-})$.
\end{theorem}
\begin{example}[Pseudosphere]
 \label{ex:pseudosphere}
 Define $f\colon{}\R^2\to\R^3$ as
 \[
    f(x,y):=
    \left(
      \sech x\, \cos y, \sech x\, \sin y,
      x - \tanh x
    \right).
 \]
 If we set 
 $\nu:=(\tanh x\,\cos y, \tanh x\,\sin y, \sech x)$,
 then $\nu$ is the unit normal vector and $f$ is a front 
 whose singular set $\{x=0\}$ consists of cuspidal edges.
 The Gaussian curvature of $f$ is $-1$, and 
 the coordinate system $(u,v)$ defined as $x=u-v$, $y=u+v$ is the 
 asymptotic Chebyshev net (see Example~\ref{ex:chebyshev}) with
 $\theta=4\arctan\exp(u-v)$.

 Since $f(x,y+2\pi)=f(x,y)$, $f$ induces a smooth map  $f_1$
 from the cylinder $M^2=\R^2/\{(0,2\pi m)\,;\,m\in\Z\}$ into $\R^3$.
 The front $f_1\colon{}M^2\to \R^3$ has two ends $p_1$, $p_2$
 with growth order $a(p_j)=0$.
 Hence by Theorem~\ref{thm:Gauss-Bonnet-complete}, we have
 \[
     2\int_{\Sigma_{f_1}}\!\!\kappa_s\,ds = \Area(M^2)=8\pi.
 \]
 In fact, the singular curvature is positive.
\end{example}
\begin{example}[Kuen's surface]
 \label{ex:kuen}
 The smooth map $f\colon{}\R^2\to\R^3$ defined as 
 \[
    f(x,y)=
    \frac{1}{1+2(1+2y^2)e^{2x}+e^{4x}}
    \begin{pmatrix}
      4 e^x(1+e^{2x})(\cos y+y\sin y)\\
      4 e^x(1+e^{2x})(\sin y+y\cos y)\\
      2 + 2x(1+2y^2)e^{2x}+(x-2)e^{4x}
    \end{pmatrix}
 \]
 is called {\em Kuen's surface}, which is considered as a
 weakly complete front with the unit normal vector
 \[
    \nu(x,y)=
    \frac{1}{1+2(1+2y^2)e^{2x}+e^{4x}}
    \begin{pmatrix}
     8e^{2x} y \cos y -
     (1+2(1-2y^2)e^{2x}+e^{4x})\sin y\\
     8e^{2x} y \sin y +
     (1+2(1-2y^2)e^{2x}+e^{4x})\cos y\\
     4e^{x}(1-e^{2x})y
    \end{pmatrix},
 \]
 and has  Gaussian curvature $-1$.
 The coordinate system $(u,v)$ such that $x=u-v$ and $y=u+v$
 is the asymptotic Chebyshev net with 
 $\theta=-4\arctan\bigl(2ye^x/(1+e^{2x})\bigr)$.
 Since the singular set $\Sigma_f=\{y=0\}\cap \{y=\pm\cosh x\}$ 
 is non-compact, $f$ is not complete.
\end{example}
\begin{example}[Cones]
 \label{ex:cone}
 Define $f\colon{}\R^2\setminus\{(0,0)\}\to \R^3$ as
 \[
    f (x,y) = (\log r\cos\theta,\log r\sin\theta,a\log r)
    \qquad (x,y) = (r\cos\theta,r\sin\theta),
 \]
 where $a\neq 0$ is a constant.
 Then $f$ is a front with $\nu=(a\cos t,a\sin t,-1)/\sqrt{1+a^2}$.
 The singular set is $\Sigma_f=\{r=1\}$, which corresponds to the single
 point $(0,0,0)\in\R^3$.
 That is, all points in $\Sigma_f$ are degenerate singular points.
 The image of the singular points is a cone of angle
 $\mu=2\pi/\sqrt{1+a^2}$
 and the area growth order of the two ends are $1/\sqrt{1+a^2}$.
  Theorem~\ref{thm:Gauss-Bonnet-complete} cannot be applied to
 this example because the singularities degenerate.
 However, this example suggests that it might be natural to define
 the ``singular curvature measure'' at a cone-like singularity
 as the cone angle.
\end{example}
\section{Behavior of the Gaussian curvature} 
\label{sec:gaussian-curvature}
Firstly, we shall prove the following assertion, which says
that the shape of
singular points is very restricted when the Gaussian curvature is
bounded.
\begin{theorem}
 \label{thm:gaussian-singular}
 Let $f\colon{}M^2\to (N^3,g)$ be a front,
 $p\in M^2$ a singular point,
 and $\gamma(t)$ a singular curve consisting of non-degenerate singular
 points with $\gamma(0)=p$  defined on an open interval $I\subset\R$.
 Then the Gaussian curvature $K$ is bounded on
 a sufficiently small neighborhood of $\gamma(I)$
 if and only if the second fundamental form vanishes on 
 $\gamma(I)$.

 Moreover, 
 if the extrinsic curvature $K_{\ext}$
 {\rm(}i.e.\ the product of the principal curvatures{\rm )}
 is non-negative on $U\setminus\gamma(I)$ for a neighborhood of 
 $U$ of $p$, 
 then the singular curvature is non-positive.
 Furthermore, if $K_{\ext}$ is bounded below by a positive constant on
 $U\setminus\gamma(I)$ then the singular curvature at $p$ takes a
 strictly negative value.

 In particular, when $(N^3,g)=(\R^3,g_0)$, the singular curvature
 is non-positive if the Gaussian curvature $K$ is non-negative
 near the singular set.
\end{theorem}
\begin{proof}%
[Proof of the first part of Theorem~\ref{thm:gaussian-singular}]
 We shall now prove the first part of the theorem.
 Take an adapted coordinate system $(u,v)$ such that the singular
 point $p$ corresponds to $(0,0)$, and 
 write the second fundamental form of $f$ as
 \begin{equation}\label{eq:second-coef}
     h=L\,du^2+2\,M\,du\,dv+N\,dv^2\quad
    \left(
    \begin{array}{r@{}l}
     L&=-g(f_u,\nu_u),~ N=-g(f_v,\nu_v),\\
     M&=-g(f_v,\nu_u)=-g(f_u,\nu_v)
    \end{array}
    \right).
 \end{equation}
 Since $f_u$ and $f_v$ are linearly dependent on the $u$-axis,
 $LN-(M)^2$ vanishes on the $u$-axis as well as the area density
 function $\lambda(u,v)$.
 Then by the Malgrange preparation theorem (see \cite[page 91]{GG}), 
 there exist smooth functions
 $\phi(u,v)$, $\psi(u,v)$ such that
 \begin{equation}\label{eq:expand-margrange}
  \lambda(u,v)=v\phi(u,v)\qquad\text{and}\qquad
   LN-(M)^2=v \psi(u,v).
 \end{equation}
 Since \eqref{eq:non-degenerate-adapted}, $\lambda_v\neq 0$ holds.
 Hence $\phi(u,v)\neq 0$ on a neighborhood of the origin.

 Firstly, we consider the case $p$ is a cuspidal edge point.
 Then we can choose $(u,v)$ so that 
 $\partial/\partial v$ gives the null direction.
 Since $f_v=0$ holds on the $u$-axis, we have 
 $M=N=0$.
 By \eqref{eq:k-int-ext} and \eqref{eq:expand-margrange}, we have
 $K=c_{N^3}+\psi(u,v)/(v \phi(u,v)^2)$.
 Thus the Gaussian curvature is bounded if and only if
 \[
    L(u,0)N_v(u,0)=
    \left.\bigl(LN-(M)^2\bigr)_v\right|_{v=0}=\psi(u,0)=0
 \]
 holds on the $u$-axis. 
 To prove the assertion, it is sufficient to show that
 $N_v(0,0)\ne 0$.
 Since 
 $\lambda_v=\mu_g(f_{u},f_{vv},\nu)\neq 0$, 
 $\{f_{u},f_{vv},\nu\}$ is linearly independent.
 Here, we have
 \[
    2g(\nu_v,\nu)=g(\nu,\nu)_v=0\qquad\text{and}\qquad
      \left. g(\nu_v,f_{u})\right|_{v=0}=-M=0. 
 \]
 Thus $\nu_v=0$ if and only if $g(\nu_v,f_{vv})=0$.
 On the other hand, $\nu_v(0,0)\ne 0$ holds, since $f$ is a front
 and $f_v=0$. Thus  we have
 \begin{equation}\label{eq:Nv-neq-zero}
     N_v(0,0)=g(f_v,\nu_v)_v=g(\nu_v,f_{vv})\ne 0.
 \end{equation}
 Hence the first part of Theorem~\ref{thm:gaussian-singular}
 is proved for cuspidal edges.

 Next we consider the case that $p$ is not a cuspidal edge point.
 Under the same notation as in the previous case,
 $f_u(0,0)=0$ holds because $p$ is not a cuspidal edge.
 Then  we have  $M(0,0)=L(0,0)=0$, and thus
 the Gaussian curvature is bounded if and only if
 \[
   L_v(u,0)N(u,0)=
    \left.\bigl(LN-(M)^2\bigr)_v\right|_{v=0}=\psi(u,0)=0
 \]
 holds on the $u$-axis. 
 Thus, to prove the assertion, it is sufficient to show that
 $L_v(0,0)\ne 0$.
 Since $\lambda_v=\mu_g(f_{uv},f_{v},\nu)$ does not vanish,
 $\{f_{uv},f_{v},\nu\}$ is linearly independent.
 On the other hand, $\nu_u(0,0)\ne 0$, 
 because $f$ is a front and $f_u(0,0)=0$.
 Since $g(\nu,\nu)_v=0$ and $g(\nu_u,f_{v})=-M=0$,
 we have
 \begin{equation}
  \label{eq:Lv-neq-0}
  L_v(0,0)=g(\nu_{u},f_{u})_v=g(\nu_{u},f_{uv})\ne 0.
 \end{equation}
 Hence the first part of the theorem is proved.
\end{proof}
Before proving the second part of Theorem~\ref{thm:gaussian-singular},
we prepare the following lemma:
\begin{lemma}%
[Existence of special adapted coordinates along cuspidal edges]
\label{lem:normal-coordinate}
 Let $p$ be a cuspidal edge of a front $f\colon{}M^2\to (N^3,g)$.
 Then there exists an adapted coordinate system
 $(u,v)$ satisfying the following properties{\rm :}
 \begin{enumerate}
  \item\label{item:normal-1} 
       $g(f_u,f_u)=1$ on the $u$-axis,
  \item\label{item:normal-2}
       $f_v$ vanishes on the $u$-axis,
  \item\label{item:normal-3} 
       $\lambda_v=1$ holds on the  $u$-axis,
  \item\label{item:normal-4} 
       $g(f_{vv},f_u)$ vanishes on  the $u$-axis, and
  \item\label{item:normal-5} 
       $\{f_u,f_{vv},\nu\}$ is a positively oriented orthonormal 
       basis along the $u$-axis.
 \end{enumerate}
\end{lemma}
We shall call such a coordinate system $(u,v)$ 
a {\em special adapted coordinate system}.
\begin{proof}[Proof of Lemma~\ref{lem:normal-coordinate}]
 One can easily take an adapted coordinate system $(u,v)$ at $p$
 satisfying \ref{item:normal-1} and \ref{item:normal-2}.
 Since $\lambda_v\neq 0$ on the $u$-axis, we can choose $(u,v)$
 as $\lambda_v>0$ on the $u$-axis.
 In this case, $r:=\sqrt{\lambda_v}$ is a smooth function on a
 neighborhood of $p$.
 Now we set
 \[
    u_1=u,\qquad v_1=\sqrt{\lambda_v(u,0)}\,v.
 \]
 Then the Jacobian matrix is given by
 \[
    \frac{\partial (u_1,v_1)}{\partial (u,v)}
        = 
    \begin{pmatrix}
        1 & 0 \\ r'(u) & r(u)
    \end{pmatrix},
  \quad \text{where $r(u):=\sqrt{\lambda_v(u,0)}$.}
 \]
 Thus we have
 \[
   \left. (f_{u_1},f_{v_1})\right |_{v=0}=
   \left.
   (f_u,f_v)
   \begin{pmatrix}
    1 & 0 \\ \frac{r'(u)}{r(u)}v & \frac{1}{r(u)}
   \end{pmatrix}
   \right |_{v=0}=
   (f_u,f_v)
   \begin{pmatrix}
    1 & 0 \\ 0 & \frac{1}{r(u)}
   \end{pmatrix}.
 \]
 This implies that
 $f_{u_1}=f_u$ and  $f_{v_1}=0$ on the $u$-axis.
 Thus the new coordinates $(u_1,v_1)$ satisfy \ref{item:normal-1}
 and \ref{item:normal-2}.
 The signed area density function with respect to $(u_1,v_1)$ is given
 by $\lambda_1:=\mu_g(f_{u_1},f_{v_1},\nu)$.
 Since $f_{v_1}=0$ on the $u$-axis, we have
 \begin{equation}
  \label{eq:e1}
  (\lambda_1)_{v_1}:=
   \mu_g(f_{u_1},D_{v_1}f_{v_1},\nu).
 \end{equation}
 On the other hand, we have
 \begin{equation}
  \label{eq:e2}
    f_{v_1}=\frac{f_v}{r(u)} \quad\text{and}\quad
  D_{v_1}f_{v_1}=\frac{D_{v}f_{v_1}}{r(u)}
    =\frac{f_{vv}}{r^2}=\frac{f_{vv}}{\lambda_v}
 \end{equation}
 on the $u_1$-axis.
 By \eqref{eq:e1} and \eqref{eq:e2}, we have
 $(\lambda_1)_{v_1}=\lambda_v/\lambda_v=1$ and
 have shown that $(u_1,v_1)$ satisfies \ref{item:normal-1},
 \ref{item:normal-2} and \ref{item:normal-3}.

 Next, we set
 \[
    u_2:=u_1+v_1^2\,s(u_1), \qquad v_2:=v_1,
 \]
 where $s(u_1)$ is a smooth function in $u_1$.
 Then we have
 \[
    \frac{\partial (u_2,v_2)}{\partial (u_1,v_1)}
    =\begin{pmatrix}
        1+v_1^2s' & 2v_1 s(u_1) \\
        0 & 1
      \end{pmatrix},
 \]
 and
 \[
    \left.\frac{\partial (u_1,v_1)}{\partial (u_2,v_2)}
    \right |_{v_1=0}
    =\left. \frac{1}{1+v_2^2s'}
     \begin{pmatrix}
      1 & -2v_1 s(u_2) \\ 0 & 1+v_2^2s'
     \end{pmatrix}
     \right |_{v_2=0}
     =\begin{pmatrix} 
       1& 0 \\ 0 & 1
      \end{pmatrix}.
 \]
 Thus the new coordinates $(u_2,v_2)$ satisfy 
 \ref{item:normal-1} and \ref{item:normal-2}.
 On the other hand,
 the area density function 
 $\lambda_2:=\mu_g(f_{u_2},f_{v_2},\nu)$
 satisfies
 \[
    (\lambda_2)_{v_2}=
          \mu_g(f_{u_2},f_{v_2},\nu)_{v_2}
                =\mu_g(f_{u_2},D_{v_2}f_{v_2},\nu).
 \]
 We have on the $u_2$-axis that $f_{u_2}=f_{u_1}$
 and
 \begin{align}
  f_{v_2}&=\frac{-2 v_1s }{1+v_1^2s'}f_{u_1}+f_{v_1}
  \label{eq:l-1} \\
  g(D_{v_2}f_{v_2})&=D_{v_1}f_{v_2}=
   \frac{-2 s }{1+v_1^2s'}f_{u_1}+D_{v_1}f_{v_1}.
  \label{eq:l-2}
 \end{align}
 Thus one can easily check that $(\lambda_2)_{v_2}=1$
 on the $u$-axis.
 By \eqref{eq:l-2}, we have
 $g(f_{u_2},D_{v_2}f_{v_2})=-2s +g(f_{u_1},D_{v_1}f_{v_1})$.
 Hence, if we set
 \[
    s(u_1):=\frac{1}{2}g\bigl(f_{u_1}(u_1,0),(D_{v_1}f_{v_1})(u_1,0)\bigr),
 \]
 then the coordinate $(u_2,v_2)$ satisfies
 \ref{item:normal-1}, \ref{item:normal-2}, \ref{item:normal-3} and 
 \ref{item:normal-4}.
 Since 
 $g\bigl((f_{v_2})_{v_2},\nu\bigr)=-g(f_{v_2},\nu_{v_2})=0$,
 $f_{v_2v_2}(u_2,0)$ is perpendicular to both $\nu$ and $f_{u_2}$.
 Moreover,  we have on the $u_2$-axis
 \[
   1=(\lambda_2)_{v_2}
     =\mu_g(f_{u_2},f_{v_2},\nu)_{v_2}
     =\mu_g(f_{u_2},D_{v_2}f_{v_2},\nu)
       =g(D_{v_2}f_{v_2},f_{u_2}\times_g \nu),
 \]
 and can conclude that $D_{v_2}f_{v_2}$ is a unit vector.
 Thus $(u_2,v_2)$ satisfies \ref{item:normal-5}.
\end{proof}
Using the existence of the special adapted
coordinate system, we shall show the
second part of the theorem.
\begin{proof}%
[Proof of the second part of Theorem~\ref{thm:gaussian-singular}]
 We suppose $K\geq c_{N^3}$, where $c_{N^3}$ is the sectional 
 curvature of $(N^3,g)$ with respect to the tangent plane.
 Then by \eqref{eq:k-int-ext}, $K_{\ext}\geq 0$ holds.

 If a given non-degenerate singular point $p$ is not
 a cuspidal edge, the singular curvature is negative by
 Corollary~\ref{cor:singular-curvature-peak}.
 Hence it is sufficient to consider the case that $p$ is a cuspidal 
 edge.
 So we may take a special adapted coordinate system as in 
 Lemma~\ref{lem:normal-coordinate}.
 We take smooth functions $\varphi$ and $\psi$ as in
 \eqref{eq:expand-margrange}.

 Since $K$ is bounded, 
 $\psi(u,0)=0$ holds, as seen in the proof of the first part.
 By the Malgrange preparation theorem again,
 we may put
 $LN-M^2=v^2\psi_1(u,v)$,
 and have the expression
 $K_{\ext}={\psi_1}/{\phi^2}$.
 Since $K_{\ext}\geq 0$, we have $\psi_1(u,0)\geq 0$.
 Moreover, if $K_{\ext}\geq \delta>0$ on a neighborhood
 of $p$, then $\psi_1(u,0)>0$.
 Since $L=M=N=0$  on the $u$-axis,  we have
 \begin{equation}\label{eq:uv-positive}
    0\leq 2\psi_1(u,0)=\bigl(LN-(M)^2\bigr)_{vv}=L_vN_v-(M_v)^2
     \leq L_vN_v.
 \end{equation}
 Here,  $\{f_u, f_{vv},\nu\}$ is an orthonormal basis,
 and  $g(f_{uu},f_u)=0$ and $L=g(f_{vv},\nu)=0$  on the $u$-axis.
 Hence
 \[
   f_{uu}=g(f_{uu},f_{vv})f_{vv}+g(f_{uu},\nu) \nu=g(f_{uu},f_{vv})f_{vv}.
 \]
 Similarly, since
 $2g(\nu_v,\nu)=g(\nu,\nu)_v=0$ and $g(\nu_v,f_u)=-M=0$, we have
 \[
    \nu_v=g(\nu_v,f_{vv}) f_{vv}.
 \]
 Since $\lambda_v=1>0$ and  $|f_u|=1$,  the singular curvature is given by
 \begin{equation}\label{eq:singular-curvature-normal}
    \kappa_s=\mu_g(f_u,f_{uu},\nu)
        =g(f_{uu},f_{vv})\mu_g(f_u,f_{vv},\nu)=g(f_{uu},f_{vv})
              =\frac{g(f_{uu},\nu_v)}{g(f_{vv},\nu_v)}.
 \end{equation}
 On the other hand, we have on the $u$-axis that
 \[
   -L_v=g(f_u,\nu_u)_v=g(f_{uv},\nu_u)+g(f_u,\nu_{uv})=
            g(f_u,\nu_{uv}),
 \]
 because
 $g(f_{uv},\nu_u)|_{v=0}=-M_u(u,0)=0$.
 Moreover, we have
 \[
   \nu_{uv}=D_vD_u \nu =
    D_uD_v \nu+R(f_v,f_u)\nu =D_uD_v \nu=\nu_{vu}
 \]
 since $f_v=0$, where $R$ is the Riemannian curvature tensor of
 $(N^3,g)$.
 Thus, 
 \[
   L_v=
   -g(f_u,\nu_{uv})
   =-g(f_u,\nu_v)_u+g(f_{uu},\nu_v)
   =M_u+g(f_{uu},\nu_v)=g(f_{uu},\nu_v)
 \]
 holds.
 Since we have on the $u$-axis that
 \[
    -N_v=g(f_v,\nu_v)_v=g(f_{vv},\nu_v)
      +g(f_v,\nu_{vv})=g(f_{vv},\nu_v),
 \]
 \eqref{eq:singular-curvature-normal} and 
 \eqref{eq:uv-positive} imply that
 \[
     \kappa_s=-\frac{L_v}{N_v}=-\frac{L_vN_v}{N_v^2}
              \leq 0.
 \]
 If $K_{\ext}\geq \delta>0$, \eqref{eq:uv-positive}
 becomes $0<L_vN_v$, and we have $\kappa_s<0$.
\end{proof}
\begin{remark}
 \label{rem:alexsandrov}
 Let $f\colon{}M^2\to \R^3$ be a compact front with positive Gaussian 
 curvature. For example, parallel surfaces of
 compact immersed constant mean curvature surfaces 
 (e.g.\ Wente tori) give such examples. 
 In this case, we have 
 the following opposite of the Cohn-Vossen inequality
 by Theorem~\ref{thm:Gauss-Bonnet-compact}:
 \[
   \int_{M^2}K dA> 2\pi\chi(M^2).
 \]
 On the other hand, the total curvature of a compact $2$-dimensional 
 Alexandrov space is bounded from above  by $2\pi\chi(M^2)$
 (see Machigashira \cite{Mac}).
 This implies that a front with positive curvature
 cannot be a limit of Riemannian $2$-manifolds with
 Gaussian curvature bounded below by a constant.
 We can give another explanation of this phenomenon as follows: 
 Since $K>0$, we have $\kappa_s<0$ and the shape of the surfaces 
 looks like cuspidal  hyperbolic parabola.
 So if the front is a limit of the
 sequence of immersions $f_n$, the curvature of $f_n$ must 
 converge to $-\infty$.
\end{remark}
\begin{example}[Fronts of constant positive Gaussian curvature]
\label{ex:positive-curvature}
 Let $f_0\colon{}M^2\to \R^3$
 be an immersion of constant mean curvature $1$
 and $\nu$ the unit normal vector of $f_0$.
 Then the parallel surface $f:=f_0-\nu$ gives a front
 of constant Gaussian curvature $1$.
 If we take isothermal principal curvature coordinates $(u,v)$ 
 on $M^2$ with respect to $f_0$, the first and
 second fundamental forms of $f$ are given by
 \[
    ds^2=dz^2+2\cosh \theta\, dz d\bar z +d\bar z^2 \qquad
    h=2\sinh\theta\, dz d \bar z,
 \]
 where $z=u+iv$ and $\theta$ is a real-valued function in $(u,v)$, 
 which is called the {\em complex Chebyshev net}.
 The sinh-Gordon equation $\theta_{uu}+\theta_{vv}+4\sinh \theta=0$
 is the integrability condition.
 In this case, the singular curve is characterized by
 $\theta=0$, and the condition for non-degenerate singular points 
 is given by $d\theta\ne 0$.
 Moreover, the cuspidal edges are characterized by
 $\theta_v\ne 0$, and the swallowtails are
 characterized by $\theta_u\ne 0, \theta_v=0$ and
 $\theta_{vv}\ne 0$.
 The singular curvature on cuspidal edges is given by
 \[
    \kappa_s=-\frac{(\theta_u)^2+(\theta_v)^2}{4 |\theta_v|}< 0.
 \]
 The negativity of $\kappa_s$ has been shown in
 Theorem~\ref{thm:gaussian-singular}.
 Like the case of fronts of constant negative curvature,
 Ishikawa-Machida \cite{IM} also showed that
 the generic singularities of fronts of constant positive 
 Gaussian curvature are cuspidal edges or
 swallowtails.
\end{example}

Here we should like to
remark on the behavior of mean curvature function
near the non-degenerate singular points.
\begin{corollary}
 Let $f:M^2\to (N^3,g)$ be a front and $p\in M^2$
 a non-degenerate singular point.
 Then the mean curvature function of $f$ is unbounded near
 $p$.
\end{corollary}
\begin{proof}
 The mean curvature function  $H$ is given by
 \[
    2H:=\frac{EN-2FM+GL}{EG-F^2}=\frac{EN-2FM+GL}{2\lambda^2}.
 \]
 We may assume that $u$-axis is a singular curve.
 By applying L'Hospital's rule, we have
 \[
   \lim_{v\to 0}H=
   \lim_{v\to 0}\frac{E_vN+EN_v-2F_vM-2FM_v+G_vL-GL_v}{2\lambda \lambda_v}.
 \]
 Firstly, we consider the case $(0,0)$ is a cuspidal edge.
 Then by the proof of the first part of
 Theorem~\ref{thm:gaussian-singular}, we have
 \[
    F(0,0)=G(0,0)=M(0,0)=N(0,0)=G_v(0,0)=0. 
 \]
 Thus
 \[
    \lim_{v\to 0}H=
    \lim_{v\to 0}\frac{EN_v}{2\lambda \lambda_v}.
 \]
 Since $\lambda(0,0)=0$ and 
 $N_v(0,0)\ne 0$ as shown in
 the proof of Theorem~\ref{thm:gaussian-singular},
 $H$ diverges. 

 Next, we consider the case that $(0,0)$ is not a cuspidal edge.
 When $p$ is not a cuspidal edge,
 by the proof of the first part of Theorem~\ref{thm:gaussian-singular}, 
 we then have
 \[
   E(0,0)=F(0,0)=L(0,0)=M(0,0)=E_v(0,0)=0, \qquad
   L_v(0,0)\ne 0.
 \]
 Thus
 \[
    \lim_{v\to 0}H=-
     \lim_{v\to 0}{GL_v}/{(2\lambda \lambda_v)}
 \]
 diverges, since
 $\lambda(0,0)=0$ and $L_v(0,0)\ne 0$.
\end{proof}

\subsection*{Generic behavior of the curvature near cuspidal edges}
As an application of Theorem~\ref{thm:gaussian-singular},
we shall investigate the generic behavior of the Gaussian curvature
near cuspidal edges and swallowtails in $(\R^3,g_0)$.

\begin{figure}
\footnotesize
\begin{center}
\begin{tabular}{c@{\hspace{7mm}}c@{\hspace{7mm}}c}
\raisebox{5mm}{
 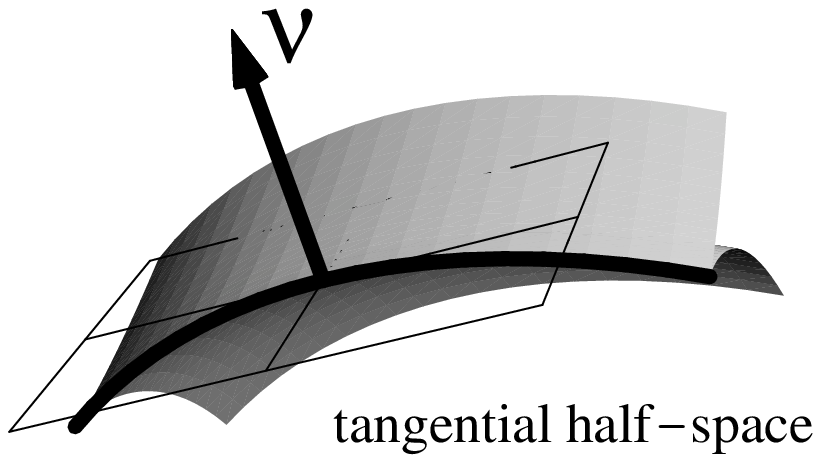}&
 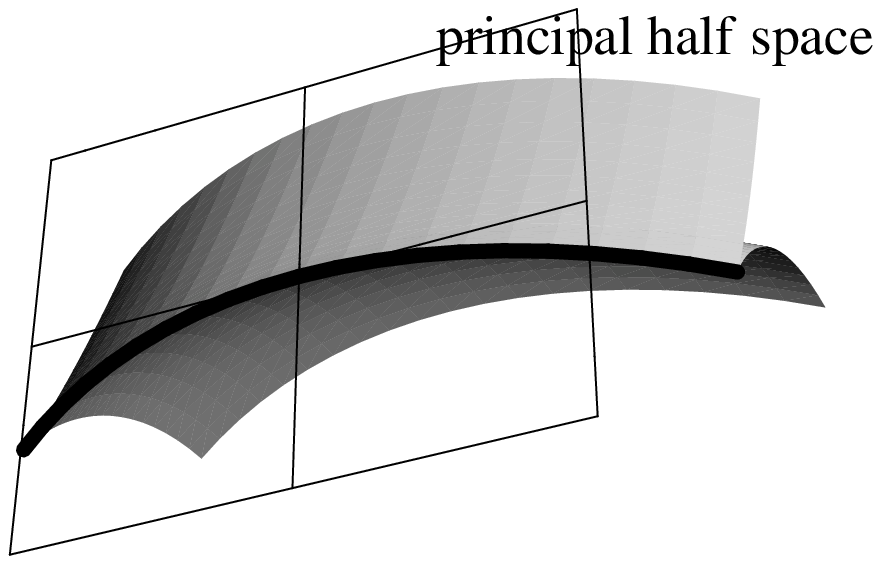 &
 \input{figure-3-1-c.pstex_t} \\
 \begin{minipage}{3.5cm}
 the half-space containing the singular curve is $\Omega(-\nu)$
 \end{minipage} &
 The outward normal & 
 the vector $\tau_0$
\end{tabular}
\end{center}
\caption{%
The half-space containing the singular curve, Theorem~\ref{thm:generic-c}.}
\label{fig:outward-normal}
\end{figure}

We call a given cuspidal edge $p\in M^2$
 {\em generic\/} if the second fundamental form does not vanish at
$p$.
Theorem~\ref{thm:gaussian-singular} implies that fronts with bounded
Gaussian curvature have only non-generic cuspidal edges.
In the proof of the theorem for cuspidal edges,
$L=0$ if and only if $f_{uu}$ is perpendicular to both $\nu$ and $f_u$,
which implies that the osculating plane of the singular curve 
coincides with the limiting tangent plane,
and we get the following:
\begin{corollary}\label{cor;cusp}
 Let $f\colon{}M^2\to \R^3$ be a front.
 Then a cuspidal edge $p\in M^2$ is 
 generic if and only if
 the osculating plane of the singular curve does not
 coincide with the limiting tangent plane at $p$.
 Moreover, the Gaussian curvature is unbounded 
and changes sign  between the two  sides of 
a generic cuspidal edge.
\end{corollary}
\begin{proof}
 By \eqref{eq:k-int-ext} and \eqref{eq:expand-margrange},
 $K=\psi/\bigl(v\varphi^2\bigr)$,
 where $\psi(0,0)\neq 0$ if $(0,0)$ is generic.
 Hence $K$ is unbounded and changes 
 sign between the two sides along the  generic cuspidal edge.
\end{proof}
We shall now determine which side has positive Gaussian curvature:
Let $\gamma$ be a singular curve of $f$ consisting of cuspidal 
edge points, and let $\hat\gamma=f\circ\gamma$.
Define
\begin{equation}\label{eq:normal-curvature}
   \kappa_{\nu}:=\frac{g_0(\hat \gamma'',\nu)}{|\hat\gamma'|^2}
\end{equation}
on the singular curve,
which is independent of the choice of parameter $t$.
We call it the {\em limiting normal curvature\/} of the cuspidal
edge $\gamma(t)$.
Then one can easily check that
{\em $p$ is a generic cuspidal edge if and only if $\kappa_{\nu}(p)$
does not vanish}.
Let $\Omega({\nu})$ (resp. $\Omega({-\nu})$) be the 
half-space bounded by the limiting tangent plane 
such that $\nu$ (resp. $-\nu$) points into 
$\Omega(\nu)$ (resp.\ $\Omega(-\nu)$).
Then the singular curve lies in $\Omega(\nu)$ if 
$\kappa_{\nu}(p)>0$ and lies in 
$\Omega(-\nu)$ if 
$\kappa_{\nu}(p)<0$.
We call $\Omega(\nu)$ (resp.\ $\Omega(-\nu)$) 
the {\em half-space containing the singular curve\/}
at the cuspidal edge point $p$.
This half-space is in general different from the principal
half-space
(see Definition~\ref{def:principal} and Figure~\ref{fig:outward-normal}).

We set
\begin{multline*}
  \sign_0(\nu):=\sign(\kappa_\nu)\\
 =
 \begin{cases}
  \hphantom{-}1 & 
  \text{
    (if $\Omega(\nu)$ is the half-space containing the singular curve)}\\
  -1 & 
  \text{
   (if $\Omega(-\nu)$ is the half-space containing the singular curve)}.
 \end{cases}
\end{multline*}
On the other hand, one can choose 
the {\em outward normal vector} {$\nu_0$} near a given cuspidal edge $p$
as in the middle figure of Figure~\ref{fig:outward-normal}.
Let $\Delta$ be a sufficiently small domain 
consisting of regular points sufficiently close to $p$ 
that lies only to one side of the cuspidal edge. 
For a given unit  normal vector $\nu$ of the front, we define
its sign $\sign_{\Delta}(\nu)$  by
$\sign_{\Delta}(\nu)=1$ (resp.\ $\sign_{\Delta}(\nu)=-1$)
if $\nu$ coincides with the outward normal $\nu_0$ on 
$\Delta$. 
The following assertion holds: 
\begin{theorem}\label{thm:generic-c}
 Let $f\colon{}M^2\to (\R^3,g_0)$ be a front, 
 $p$ a cuspidal edge and
 $\Delta$ a sufficiently small domain consisting of regular points
 sufficiently close to $p$ that lies only to one side of the cuspidal
 edge. 
 Then $\sign_{\Delta}(\nu)$ coincides with the sign of 
 the function 
 $g_0({\hat\sigma}'',{\hat \nu}')$
 at $p$, namely
\begin{equation}
  \label{eq:domain-sign}
 \sign_{\Delta}(\nu)= \sign g_0({\hat\sigma}'',{\hat \nu}')
  \qquad \left({}'=\frac{d}{ds}, ~''=D_s\frac{d}{ds}\right),
\end{equation}
 where $\sigma(s)$ is an arbitrarily fixed null curve starting 
 at $p$ and moving into  $\Delta$, and
 $\hat\sigma(s)=f(\sigma(s))$ and $\hat \nu=\nu(\sigma(s))$.
 Moreover, if $p$ is a generic cuspidal edge, then
 \[
   \sign_0({\nu})\cdot \sign_{\Delta}(\nu)
 \]
 coincides with the sign of the Gaussian curvature on $\Delta$.
\end{theorem}
\begin{proof}
 We take a special adapted coordinate system $(u,v)$ as in
 Lemma~\ref{lem:normal-coordinate} at the cuspidal edge.
 The vector  $\tau_0:=-f_{vv}=f_u\times \nu$ lies in the limiting 
 tangent plane and points in the opposite direction of the image of the
 null  curve (see Figure~\ref{fig:outward-normal}, 
 right side).

 Without loss of generality, we may assume that $\Delta=\{v>0\}$.
 The unit normal $\nu$ is the outward normal on $\Delta$
 if and only if  $g_0(\nu_{v},\tau_0)>0$, namely
 $N_v=-g_0(f_{vv},\nu_v)>0$.
 Thus we have
 $\sign_{\{v>0\}}(\nu)=\sign(N_v)$,
 which proves \eqref{eq:domain-sign}.
 Since $p$ is generic, we have $\kappa_{\nu}(p)\ne 0$ and
 $\kappa_{\nu}(p)=L$  holds. 
 On the other hand, the sign of $K$ on $v>0$
 is equal to the sign of
 \[
   \left.\bigl(LN-(M)^2\bigr)_v \right |_{v=0}=L(u,0)N_v(u,0), 
 \]
 which proves the assertion.
\end{proof}
\begin{example}
\label{ex:parabola2}
 Consider again the cuspidal parabola 
 $f(u,v)$ as in Example~\ref{ex:parabola}.
 Then $(u,v)$ gives an adapted coordinate system
 so that $\partial/\partial v$ gives a null direction, and
 we have 
 \[
    L = g_0(f_{uu},\nu) = \frac{-2ab}{\sqrt{1+b^2(1+4a^2u^2)}},\qquad
    N_v=\frac{6}{\sqrt{1+b^2(1+4a^2u^2)}}>0.
 \]
 The cuspidal edges are generic if and only if $ab\neq 0$.
 In this case, let $\Delta$ be a domain in the upper half-plane
 $\{(u,v)\,;\,v>0\}$.
 Then the unit normal vector \eqref{eq:parabola-normal} is
 the outward normal to the cuspidal edge, that is,
 $\sign_{\Delta}(\nu)=+1$.
 The limiting normal curvature as in \eqref{eq:normal-curvature} is
 computed as $\kappa_\nu = -ab/(2|a|^2\sqrt{1+b^2(1+4a^2u^2)})$,
 and hence $\sign_0(\nu)=-\sign(ab)$.
 Then $\sign(K)=-\sign(ab)$ holds on the upper half-plane.
 In fact, the Gaussian curvature is computed as
 \[
    K = \frac{-12(ab+3av)}{v(4+\bigl(1+4a^2u^2)(4b^2+12bv+9v^2)\bigr)}.
 \]

 On the other hand,  the Gaussian curvature is bounded if $b=0$.
 Moreover, the Gaussian curvature is positive if  $a<0$.
 In this case the singular curvature
 is negative when $a<0$,
 as stated in Theorem~\ref{thm:gaussian-singular}.
\end{example}
\subsection*{Generic behavior of the curvature near swallowtails}
\begin{figure}
\footnotesize
 \begin{center}
  \footnotesize
   \begin{tabular}{c@{\hspace{2cm}}c}
        \includegraphics[width=3.8cm]{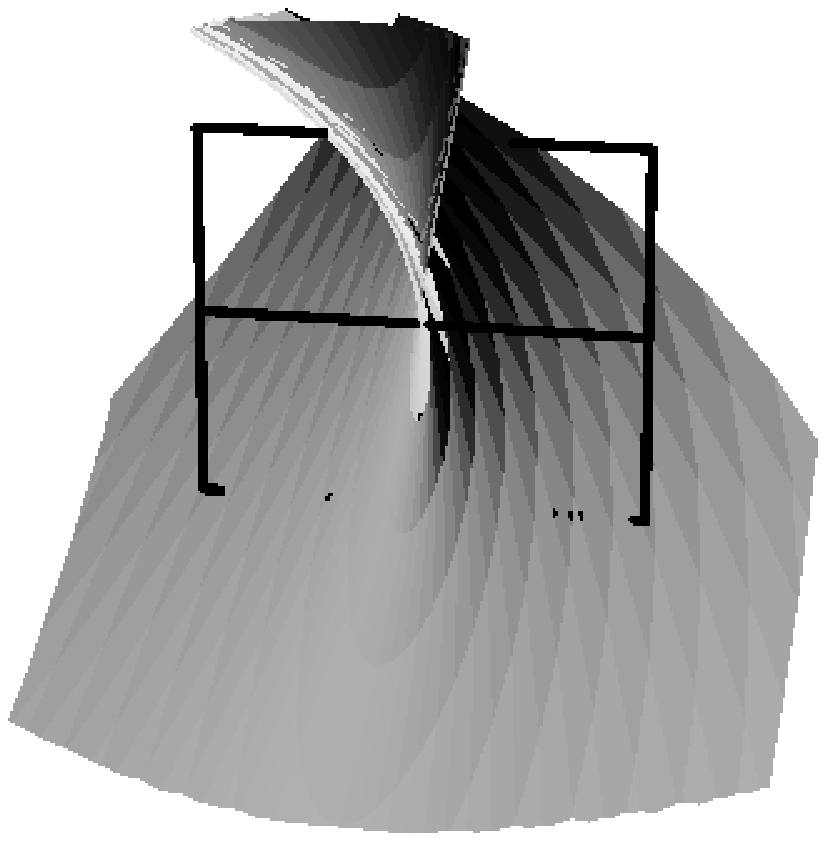} &
    \raisebox{-5mm}{
        \includegraphics[width=3.8cm]{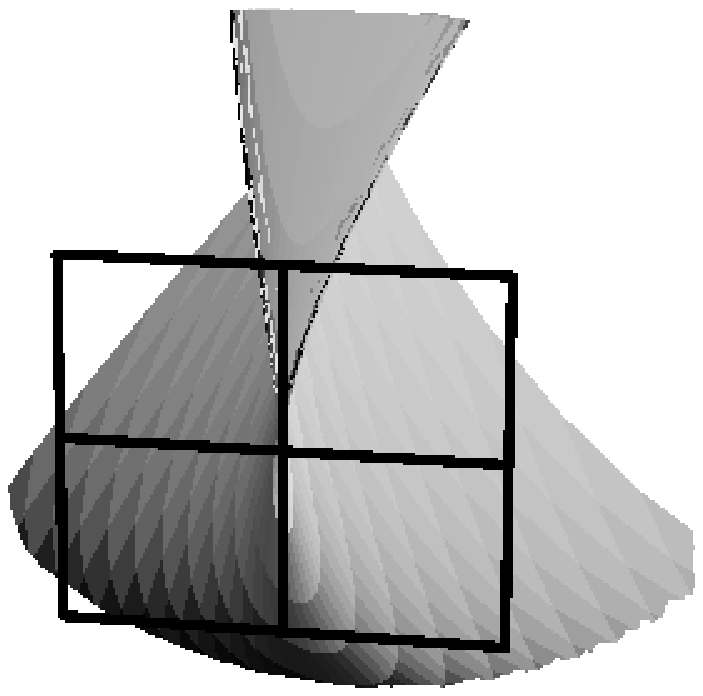}}\\
    the swallowtail $f_+$ &
    the swallowtail $f_-$ \\
   \end{tabular}
 \end{center}
 The half-space containing the singular curve is 
 the closer side 
 of the limiting tangent  plane for the left-hand figure, 
 and the farther side for the right-hand figure.
\caption{The half-space containing the singular curve for generic swallowtails
 (Example~\ref{ex:positive-negative-swallowtail}).}
\label{fig:generic-swallowtails}
\end{figure}

We call a given swallowtail $p\in M^2$ of a front 
$f\colon{}M^2\to (\R^3,g_0)$
 {\it generic\/} if the second fundamental form does not vanish at
$p$.
\begin{proposition}
 \label{prop:swallow-half-plane}
 Let $f\colon{}M^2\to (\R^3,g_0)$
 be a front and $p$ a generic swallowtail
 Then we can take a half-space $H\subset\R^3$
 bounded by the limiting tangent plane such that
 any null curve at $p$ lies in $H$ near $p$
 {\rm (}see Figure~\ref{fig:generic-swallowtails}{\rm)}.
\end{proposition}
We shall call $H$  the {\em half-space containing the singular curve\/}
at the generic swallowtail.
At the end of this section, we shall see that 
the singular curve is in fact contained in this half-space
for a neighborhood of the swallowtail 
(see Figure~\ref{fig:generic-swallowtails} and
Corollary~\ref{cor:limit-tangent}).
For a given unit  normal vector $\nu$ of the front, 
we define the sign
$\sign_0(\nu)$ of it by
$\sign_0(\nu)=1$ (resp.\ $\sign_0(\nu)=-1$)
if $\nu$ points (resp.\ does not point)
into the half-space containing the singular curve. 
\begin{proof}[Proof of Proposition~\ref{prop:swallow-half-plane}]
 Take an adapted coordinate system $(u,v)$ and assume $f(0,0)=0$ by
 translating in $\R^3$ if necessary.
 Write the second fundamental form as in \eqref{eq:second-coef}.
 Since $f_u(0,0)=0$, we have $L(0,0)=M(0,0)=0$, and 
 we have the following Taylor expansion:
 \[
    g_0\bigl(f(u,v),\nu\bigr)=
         \frac{v^2}{2} g_0\bigl(f_{vv}(0,0),\nu(0,0)\bigr) +o(u^2+v^2)
      =  \frac{1}{2}N(0,0) v^2  +o(u^2+v^2).
 \]
 Thus the assertion holds.
 Moreover we have
 \begin{equation}\label{eq:primary-signature}
    \sign(N)=\sign_0(\nu).
 \end{equation}
\end{proof}
\begin{corollary}
 \label{cor:primary-signature}
 Let $\sigma(s)$ be an arbitrary curve starting at 
 the swallowtail such that $\sigma'(0)$ is transversal to
 the singular direction.
 Then 
 \[
   \sign_0(\nu)=\sign\bigl(g_0({\hat\sigma}''(0),\nu(0,0))\bigr)
 \]
 holds¡¤where $\hat\sigma=f\circ\sigma$.
\end{corollary}

We let $\Delta$ be a sufficiently small domain consisting
of regular points sufficiently close to a swallowtail $p$. 
The domain $\Delta$ is called the {\em tail part\/} 
if $\Delta$ is on the opposite side of the self-intersection of 
the swallowtail.
We define $\sign_{\Delta}(\nu)$ by $\sign_{\Delta}(\nu)=1$
(resp.\ $\sign_{\Delta}(\nu)=-1$)
if $\nu$ is (resp.\ is not) the outward normal of $\Delta$. 
Now we have the following assertion:
\begin{theorem}
 \label{thm:secondary-signature}
 Let $f\colon{}M^2\to (\R^3,g_0)$ be a front, 
 $p$ a generic swallowtail and $\Delta$ a sufficiently small
 domain consisting of regular points sufficiently close to $p$.
 Then the Gaussian curvature is unbounded and changes  
 sign between the two sides along the singular curve.
 Moreover,
 $\sign_0(\nu)\sign_{\Delta}(\nu)$
 coincides with the sign of the Gaussian curvature
 on $\Delta$.
\end{theorem}
\begin{proof}
 If we change $\Delta$ to the opposite side, 
 $\sign_{\Delta}(\nu)\sign_{\Delta}(K)$ does not change sign.
 So we may assume that $\Delta$ is the tail part.
 We take an adapted coordinate system $(u,v)$ at the swallowtail
 and write the null vector field as
 $\eta(u)=(\partial/\partial u)+e(u)(\partial/\partial v)$,
 where $e(u)$ is a smooth function.
 Then 
 \[
     f_u(u,0) + e(u)f_v(u,0)=0\quad
     \text{and}\quad
     f_{uu}(u,0)+e_u(u)f_v(u,0)+e(u)f_{uv}(u,0)=0
 \]
 hold.
 Since $u=0$ is a swallowtail, $e(0)=0$ and $e'(0)\neq 0$ hold,
 where $'=d/du$.

 The vector $f_{uu}$ points toward the tail part $\Delta$.
 Thus $f_v$ points toward $\Delta$ if and only if 
 $g_0(f_v,f_{uu})$ is positive.
 Since $f_u=-e(u)f_v$ and $e(0)=0$, we have
 $f_{uu}(0,0)=e'(0) f_v(0,0)$ and 
 \[
   g_0\bigl(f_{uu}(0,0),f_v(0,0)\bigr)=
     -e'(0)\, g_0\bigl( f_v(0,0),f_v(0,0)\bigr).
 \]
 Thus  $g_0\bigl(f_{uu}(0,0),f_v(0,0)\bigr)$ is positive (that
 is,
 the tail part is $v>0$)
 if and only if  $e'(0)<0$.

 Changing $v$ to $-v$ if necessary, we assume $e'(0)>0$,
 that is,  the tail part lies in $v>0$.
 For each fixed value of $u\neq 0$, we take a curve 
 \[
  \sigma(s)=\bigl(u+\varepsilon  s, s|e(u)|\bigr)
             =\bigl(u+\varepsilon s, \varepsilon e(u)\,s\bigr)
   \qquad \varepsilon = \sign e(u)
 \]
 and let $\hat\sigma= f\circ\sigma$.
 Then $\sigma$ is traveling into the upper half-plane $\{v>0\}$,
 that is, $\hat\sigma$ is traveling into $\Delta$.
 Here, we have
 \begin{align*}
   {\hat\sigma}'(0)
   &= \varepsilon\bigl(f_u(u,0) + e(u)f_v(u,0)\bigr) = 0 
   \qquad\text{and}\\
   {\hat\sigma}''(0) &=
     \left.\varepsilon\bigl( \varepsilon (f_u+ef_v)_u + 
       \varepsilon e (f_u+ef_v)_v\bigr)\right|_{v=0}\\
  &= e(u) \bigl(f_{uv}(u,0)+e(u)f_{vv}(u,0)),
 \end{align*}
 where ${}'{~}=d/ds$.
 In particular, $\sigma$ is a null curve 
  starting at $(u,0)$ and traveling into $\Delta$.
 Then by Theorem~\ref{thm:generic-c}, we have
 \[
   \sign_{\Delta}(\nu)=\lim_{u\to 0}
       \sign\bigl(g_0({\hat\sigma}''_u(s),{\hat\nu}'(s))\bigr).
 \]
 Here, the derivative of $\hat\nu(t)=\nu(\sigma(t))$ is
 computed as
 ${\hat\nu}'=\varepsilon\{\nu_u(u,0)+e(u)\nu_v(u,0)\}$.
 Since $e(0)=0$, we have 
 \[\left.
   g\bigl({\hat\sigma}''(s),{\hat \nu}'(s)\bigr)
  \right|_{s=0}=
      |e(u)| g_0\bigl(f_{uv}(u,0), \nu_u(u,0)\bigr)
      +\{e(u)\}^2\varphi(u),
 \]
 where $\varphi(u)$ is a smooth function in $u$.
 Then we have
 \[
   \sign_{\Delta}(\nu)=\lim_{u\to 0}
         \sign\bigl(g_0({\hat\sigma}''(s),
                      {\hat \nu}'(s))\bigr)
         =\sign\bigl(g_0(f_{uv}(0,0), \nu_u(0,0)\bigr).
 \]
 Here, 
 $L_v(0,0)=-g_0(f_u,\nu_u)_v =-g_0(f_{uv},\nu_u)$ because $f_u=0$,
 which implies that
 \[
  \sign_{\Delta}(\nu)=\sign(L_v(0,0)).
 \]
 On the other hand, the sign of $K$ on $v>0$
 is equal to the sign of
 \[
   \left.\bigl(LN-(M)^2\bigr)_v \right |_{v=0}=
  N(0,0)L_v(0,0). 
 \]
 Then \eqref{eq:primary-signature} implies the assertion.
\end{proof}
\begin{example}
 \label{ex:positive-negative-swallowtail}
 Let
 \[
   f_{\pm}(u,v)=\frac{1}{12}
     (3u^4-12u^2v\pm (6u^2-12v)^2,8u^3-24uv,6u^2-12v).
 \]
 Then one can see that $f_{\pm}$ is a front and $(0,0)$
 is a swallowtail 
 with the unit normal vector
 \begin{multline*}
    \nu_{\pm}=\frac{1}{\delta}
    \bigl(1,u,u^2\pm12(2v-u^2)\bigr )\\
   \left(\delta =  
         \sqrt{1+u^2+145u^4+576v(v-u^2)\pm 24u^2(2v-u^2)}\right).
 \end{multline*}
 In particular, $(u,v)$ is an adapted coordinate system. 
 Since the second fundamental form is  $\pm 24 dv^2$ at the origin,  the
 swallowtail is generic and 
 $\sign_0(\nu_{\pm})=\pm 1$ because of \eqref{eq:primary-signature}.
 The images of $f_{\pm}$ are shown in Figure~\ref{fig:generic-swallowtails}.
 Moreover, since $L_v=\pm 2$ at the origin,
 $\sign_D(\nu_{\pm})=\pm 1$.
 Then by Theorem~\ref{thm:secondary-signature},
 the Gaussian curvature of the tail side of $f_{+}$  
 (resp.\ $f_{-}$) is positive (resp.\ negative).
\end{example}

Summing up the previous two theorems, we get the following:

\begin{corollary}\label{cor:limit-tangent}
 Let $\gamma(t)$ be a singular curve such that $\gamma(0)$ is a
 swallowtail.
 Then the half-space containing the singular curve at $\gamma(t)$
 converges to the half-space
 at the swallowtail $\gamma(0)$ as $t\to 0$.
\end{corollary}

\section{Zigzag numbers}
\label{sec:zigzag}
In this section, we introduce a geometric formula for
a topological invariant called the {\em zigzag number\/}.
We remark that 
Langevin, Levitt and Rosenberg \cite{LLR}
gave topological upper bounds of zig-zag numbers 
for generic compact fronts in $\R^3$.
(See Remark \ref{rem:add}.)
\subsection*{Zigzag number for fronts in the plane}
\begin{figure}
\footnotesize
\begin{center}
 \begin{tabular}{c@{\hspace{1cm}}c@{\hspace{5mm}}c}
 \input{figure-4-1-a.pstex_t} &
 \raisebox{5mm}{
  \includegraphics[height=1.5cm]{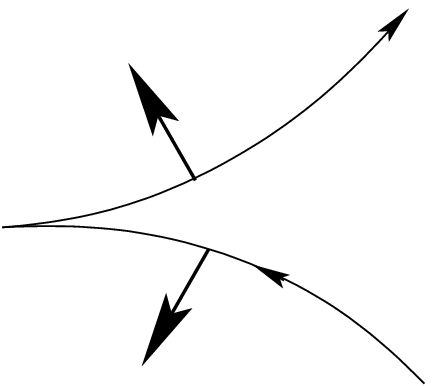}} &
 \raisebox{5mm}{
  \includegraphics[height=1.5cm]{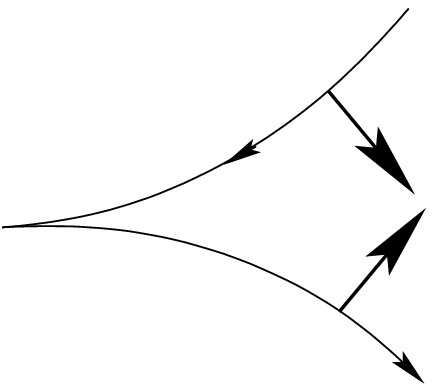}}\\
    & zig & zag \\
  (a) Canonical forms of plane fronts &
  \multicolumn{2}{c}{(b) Zig and zag for plane fronts}
 \end{tabular}\\[3mm]
 \begin{tabular}{c@{\hspace{2cm}}cc}
  \input{figure-4-1-c.pstex_t} & 
  \raisebox{5mm}{
  \includegraphics[height=2cm]{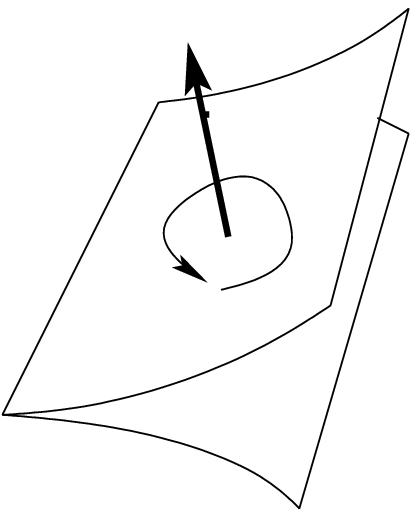}}&
  \raisebox{5mm}{
  \includegraphics[height=2cm]{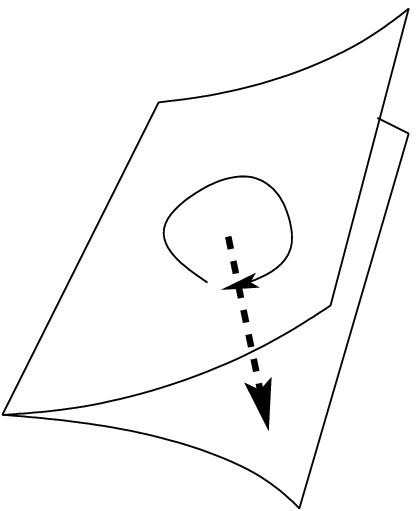}}\\
    & zig  & zag \\
  (c) Counterclockwise in $P^1(\R)$ &
  \multicolumn{2}{c}{(d) Zig and zag for cuspidal edges}
 \end{tabular}
\end{center}
\caption{Zigzag.}
\label{fig:zig-zag}
\end{figure}

First, we mention the Maslov index (see \cite{A};
which is also called the zigzag number) for fronts in the Euclidean
plane $(\R^2,g_0)$.
Let $\gamma\colon{}S^1\to\R^2$ be a generic front, that is,
all self-intersections and singularities are double points
and $3/2$-cusps, and let $\nu$ be the unit normal vector
field of $\gamma$.
Then $\gamma$ is Legendrian isotropic (isotropic as the Legendrian
lift $(\gamma,\nu)\colon{}S^1\to T_1\R^2\simeq \R^2\times S^1$) to 
one of the fronts in Figure~\ref{fig:zig-zag} (a).
The non-negative integer $m$ is called the {\em rotation number},
which is the rotational index of the unit normal vector field 
$\nu\colon{}S^1\to S^1$.
The number $k$ is called the {\em Maslov index\/}
or {\em zigzag number}.
We shall give a precise definition and a formula to calculate the 
number:
a $3/2$-cusp  $\gamma(t_0)$ of $\gamma$ is called {\em zig\/} 
(resp.\ {\em zag}) if the leftward normal vector of $\gamma$ 
points to the outside (resp.\ inside) of the cusp (see
Figure~\ref{fig:zig-zag} (b)).
We define a $C^\infty$-function $\lambda$ on $S^1$ as 
$\lambda:=\det(\gamma',\nu)$, where $'=d/dt$.
Then the leftward normal vector is given by $(\sign\lambda) \nu_{0}$.
Since $\gamma''(t_0)$ points to the inside of the cusp, 
$t_0$ is zig (resp.\ zag) if and only if
\begin{equation}\label{eq:zigzag-criterion}
    \sign\bigl(\lambda'g_0(\gamma'',\nu')\bigr) <0 \qquad
    (\text{resp.}~>0).
\end{equation}

Let $\{t_0,t_1,\dots,t_l\}$ be the set of singular points of $\gamma$
ordered by their appearance, 
and define $\zeta_j=a$ (resp.\ $=b$) if $\gamma(t_j)$ is zig (resp.\ zag),
and set 
$\zeta_{\gamma}:=\zeta_0\zeta_1\dots \zeta_l$, which is a word 
consisting of the letters $a$ and $b$.
The projection of $\zeta_{\gamma}$ to the free product $\Z_2*\Z_2$
(reduction with the relation $a^2=b^2=1$) 
is of the form $(ab)^k$ or $(ba)^k$.
The non-negative integer $k_{\gamma}:=k$ is called 
the {\em zigzag number\/}
of $\gamma$.
We shall give a geometric formula for the zigzag number
via the curvature map defined by the second author:
\begin{definition}[\cite{U}]\label{def:curvature-map}
 Let $\gamma\colon{}S^1\to \R^2$ be a front with unit normal vector $\nu$.
 The {\em curvature map\/} of $\gamma$ is the map 
 \[
    \kappa_{\gamma}\colon{}
       S^1\setminus \Sigma_{\gamma}\ni t\longmapsto
       \left[
        g_0(\gamma',\gamma'):
        g_0(\gamma',\nu')
       \right]\in P^1(\R),
 \]
 where $'=d/dt$,
 $\Sigma_{\gamma}\subset S^1$ is the set of singular points of $\gamma$,
 and $[~:~]$ denotes the homogeneous coordinates of $P^1(\R)$.
\end{definition}
\begin{proposition}\label{prop:plane-zigzag}
 Let $\gamma$ be a generic front with unit normal vector $\nu$.
 Then the curvature map $\kappa_{\gamma}$ can be extended to a
 smooth map on $S^1$.
 Moreover, the rotation number of $\kappa_{\gamma}$ is the zigzag number
 of $\gamma$.
\end{proposition}
\begin{proof}
 Let $t_0$ be a singular point of $\gamma$.
 Since $\gamma$ is a front, $\nu'(t)\neq 0$ holds on a neighborhood of
 $t_0$.
 As $\nu'$ is perpendicular to $\nu$, we have $\det(\nu,\nu')\neq 0$.
 Here, using $\lambda=\det(\gamma',\nu)$, we have 
 $\gamma'=-(\lambda/\det(\nu,\nu'))\nu'$.
 Hence we have 
 \[
    \kappa_{\gamma}=[g_0(\gamma',\gamma'):g_0(\gamma',\nu')]
    =
    \left[\lambda^2:-\frac{\lambda g_0(\nu',\nu')}{\det(\nu,\nu')}\right]
    =\left[\lambda:-\frac{g_0(\nu',\nu')}{\det(\nu,\nu')}\right]
 \]
well-defined on a neighborhood of $t_0$.
 Moreover, $\kappa_{\gamma}(t)=[0:1](=\infty)$ if and only if $t$ 
 is a singular point.
Here, we choose an inhomogeneous coordinate of $[x:y]$ as $y/x$.
 
 Since $g_0(\gamma',\nu')' =  g_0(\gamma'',\nu')$ holds 
 at a singular point $t_0$,
 $\kappa_{\gamma}$ passes through $[0:1]$ with counterclockwise 
 (resp.\ clockwise) direction if $g_0(\gamma'',\nu')>0$
 (resp.\ $<0$), see Figure~\ref{fig:zig-zag} (c).

 Let $t_0$ and $t_1$ be two adjacent zigs, and suppose $\lambda'(t_0)>0$.
 Since $\lambda$ changes sign on each cusp, we have $\lambda'(t_1)<0$.
 Then by \eqref{eq:zig-zag-criterion}, $g_0(\gamma'',\nu')(t_0)>0$
 and $g_0(\gamma'',\nu')(t_1)<0$.
 Hence $\kappa_{\gamma}$ passes through $[0:1]$  in the  
 counterclockwise direction at $t_0$, and the
 clockwise direction at $t_1$.
 Thus, this interval does not contribute to the rotation number of
 $\kappa_{\gamma}$.
 On the other hand, if $t_0$ and $t_1$ are zig and zag respectively,
 $\kappa_{\gamma}$ passes through $[0:1]$ counterclockwisely
 at both $t_0$ and $t_1$.
 Then the rotation number of $\kappa_{\gamma}$ is $1$ 
 on the interval $[t_0,t_1]$.
 Summing up, the proposition holds.
\end{proof}

\subsection*{Zigzag number for fronts in Riemannian $3$-manifolds}
Let $M^2$ be a manifold and
$f\colon{}M^2 \to N^3$ be a front with unit normal vector $\nu$
into a Riemannian 3-manifold $(N^3,g)$.
Let
$\Sigma_f\subset M^2$  be the singular set, and 
${\nu}_0$ be the unit normal vector field of $f$ defined on
$M^2\setminus\Sigma_f$ which is compatible with the orientations of $M^2$
and $N^3$, that is,
${\nu}_0= (f_u\times_g f_v)/|f_u\times_g f_v|$, where $(u,v)$ 
is a local coordinate system on $M^2$ compatible to the orientation.
Then ${\nu}_0(p)$  is $\nu(p)$ if $p\in M_+$ and $-\nu(p)$ if $p\in M_-$.

We assume all singular points of $f$ are non-degenerate.
Then each connected component $C\subset\Sigma_f$ must be a 
regular curve on $M^2$.
Let $p\in C$ be a cuspidal edge.
Then $p$ is called {\em zig\/} (resp.\ {\em zag}) if ${\nu}_0$ points
towards the outward (resp.\ inward) side of the cuspidal edge
(see Figure~\ref{fig:zig-zag} (d)).
As this definition does not depend on $p\in C$, we call
$C$  {\it zig} (resp.\ {\it zag}) if $p\in C$ is zig (resp.\ zag).

Now, we define the zigzag number for loops on $M^2$.
Take a {\em null loop\/} $\sigma\colon{}S^1\to M^2$, that is, 
the intersection of $\sigma(S^1)$ and $\Sigma_f$ consists of cuspidal
edges
and $\sigma'$ points in the null direction  at each singular point.
We remark 
that there exists a null loop in each homotopy class.
Let $Z_{\sigma}=\{t_0,\dots, t_l\}\subset S^1$ be the set of
singular points of $\sigma$ ordered by their appearance along the loop.
Define $\zeta_j=a$ (resp.\ $b$) if $\sigma(t_j)$ is zig (resp.\ zag),
and set $\zeta_{\sigma}:=\zeta_0\zeta_1\dots \zeta_l$, which is a word 
consisting of the letters $a$ and $b$.
The projection of $\zeta_{\sigma}$ to the free product $\Z_2*\Z_2$
(reduction with the relation $a^2=b^2=1$) is of the form $(ab)^k$ or $(ba)^k$.
The non-negative integer $k_{\sigma}:=k$ is called the {\em zigzag number\/}
of $\sigma$.

It is known that the zigzag number is a homotopy invariant, 
and the greatest common divisor $k_f$ of 
$\{k_{\sigma}\,| \text{$\sigma$ is a null loop on $M^2$}\}$ is the 
{\em zigzag number of $f$} (see \cite{LLR}).

\begin{remark}%
[Langevin-Levitt-Rosenberg's inequality \cite{LLR}]
\label{rem:add}
 Let $M^2$ be a compact orientable 2-manifold of genus $g$ and 
 $f:M^2\to N^3$ a front. When $N^3=\R^3$,
 \cite{LLR} proved  the following inequality
 \begin{equation} \label{eq:add}
  a_f+\frac{q_f}{2}\ge \frac{\chi_{\E}^{}}{2}+1-g +2k_f,
 \end{equation}
 where $a_f$ is the number of the connected components of 
 the singular set $\Sigma_f,$\,\, $q_f$ the number of the 
 swallowtails, and half the Euler number of the 
 limiting tangent bundle $\chi_{\E}^{}/{2}$ is 
 equal to the degree of the Gauss map.  
 Their proof is valid for the general case
 and \eqref{eq:add} holds for any $N^3$.
\end{remark}

In this section, we shall give a geometric formula for zigzag numbers 
of loops.
First, we define the normal curvature map, similar to the curvature
map for fronts in $\R^2$:
\begin{definition}[Normal curvature map]\label{def:front-curvature-map}
 Let $f\colon{}M^2\to (N^3,g)$ be a front with unit normal vector $\nu$ 
 and $\sigma\colon{}S^1\to M^2$ a null loop.
 The {\em normal curvature map\/} of $\sigma$ is the map
 \[
    \kappa_{\sigma}\colon{}
       S^1\setminus Z_{\gamma}\ni t\longmapsto
       \left[
        g(\hat\sigma',\hat\sigma'):
        g(\hat\sigma',\hat\nu')
       \right]\in P^1(\R),
 \]
 where $\hat\sigma=f\circ \sigma$, $\hat\nu=\nu\circ \sigma$,
 $'=d/dt$,
 $Z_{\sigma}\subset S^1$ is the set of singular points  of $\sigma$,
 and $[~:~]$ denotes the homogeneous coordinates of $P^1(\R)$.
\end{definition}
Then we have the following:
\begin{theorem}[Geometric formula for zigzag numbers]
\label{thm:zigzag-front}
 Let $f\colon{}M^2\to (N^3,g)$ be a front with unit normal
 vector $\nu$,
 whose singular points are all non-degenerate,
 and $\sigma\colon{}S^1\to M^2$ a null loop.
 Then the normal curvature map $\kappa_{\sigma}$ can be extended
 to $S^1$, and the rotation number of $\kappa_{\sigma}$ is equal
 to the zigzag number of $\sigma$.
\end{theorem}
\begin{proof}
 Let $t_0$ be a singular point of $\sigma$, and take a normalized coordinate
 system $(u,v)$ of $M^2$ on a neighborhood $U$ of $\sigma(t_0)$.
 Then $f_v=0$ and $f_{vv}\neq 0$ holds on the $u$-axis, and by the
 Malgrange preparation theorem, there exists a smooth function $\alpha$
 such that $g(f_v,f_v)=v^2\alpha(u,v)$ and
 $\alpha(u,0)\neq 0$.
 On the other hand, $g(f_v,\nu_v)=-N$
 vanishes  and $N_v\neq 0$ on the $u$-axis.
 Hence there exists a function $\beta$ such that
 $g(f_v,\nu_v)=v \mu(u,v)$ and $\mu(u,0)\neq 0$.
 Thus 
 \begin{equation}\label{eq:normal-curvature-extend}
   \kappa_{\sigma}=[g(f_v,f_v):g(f_v,\nu_v))]
                   =[v^2\alpha(u,v):v\beta(u,v)] 
                   =[v\alpha(u,v):\beta(u,v)] 
 \end{equation}
 can be extended to the singular point $v=0$.
 Namely, $\kappa_{\sigma}(t_0)=[0:1](=\infty)$, where 
 we choose an inhomogeneous coordinate $y/x$ for $[x:y]$.
 Moreover, $g(\hat\sigma',\hat\sigma')\neq 0$ on regular points,
 and $\kappa_{\sigma}(t)=[0:1]$ if and only if $t$ is a singular point.

 Since ${\nu}=(\sign\lambda){\nu}_0$, 
 so a singular point $t_0$
 is zig  (resp.\ zag) if and only if
 \[
     \sign(\lambda) \sign_{\Delta}(\nu)>0\qquad (\text{resp.\ }<0),
 \]
 where  $\varepsilon$ is a sufficiently small number
and  $\Delta$ is a domain containing $\sigma(t_0+\varepsilon)$ 
which lies only to one side of the cuspidal edge.
 By Theorem~\ref{thm:generic-c}, 
 $\sign_{\Delta}(\nu)=\sign g(\hat\sigma'',\hat\nu')$,
 $t_0$ is zig (resp.\ zag) if and only if
 \begin{equation}\label{eq:zig-zag-criterion}
     \sign\bigl(\hat\lambda' g(\hat\sigma'',\hat\nu')\bigr)
     >0 \qquad (\text{resp.}\ <0),
 \end{equation}
 where $\hat\lambda=\lambda\circ\sigma$.
 Since $g(\hat\sigma',\hat\nu')=g(\hat\sigma'',\hat\nu')$ holds 
 at singular points,
 we have
 \begin{itemize}
  \item if $t_0$ is zig and $\hat\lambda'(t_0)>0$ (resp.\ $<0$),
	then $\kappa_{\sigma}$ passes through $[0:1]$ counterclockwisely
	(resp.\ clockwisely).
  \item if $t_0$ is zag and $\hat\lambda'(t_0)>0$ (resp.\ $<0$),
	then $\kappa_{\sigma}$ passes through $[0:1]$ clockwisely
	(resp.\ counterclockwisely).
 \end{itemize}	
 Let $Z_{\sigma}=\{t_0,\dots,t_l\}$ be the set of singular points.
 Since the function $\lambda$ has alternative sign on the 
 adjacent domains, $\hat\lambda'(t_j)$ and $\hat\lambda'(t_{j+1})$ 
 have opposite sign.
 Thus, 
 if both $t_j$ and $t_{j+1}$ are zigs and $\hat\lambda(t_j)>0$,
 $\kappa_{\sigma}$ passes through $[0:1]$ counterclockwisely
 (resp.\ clockwisely) at $t=t_j$ (resp.\ $t_{j+1}$).
 Hence the interval $[t_j,t_{j+1}]$ does not contribute 
 to the rotation number of $\kappa_{\sigma}$.
 Similarly, two consecutive zags do not affect the rotation number.
 On the other hand, 
 if $t_j$ is zig  and $t_{j+1}$ is zag and $\hat\lambda(t_j)>0$,
 $\kappa_{\sigma}$ passes through $[0:1]$ counterclockwisely at both
 $t_j$ and $t_{j+1}$.
 Hence the rotation number of $\kappa_{\sigma}$ on the interval
 $[t_j,t_{j+1}]$ is $1$.
 Similarly, two consecutive zags increases the rotation number by $1$.
 Hence we have the conclusion.
\end{proof}
\section{Singularities of hypersurfaces}
\label{sec:hyper}
In this section, we shall investigate the behavior of sectional
curvature on fronts that are hypersurfaces.
Let $U^n$ ($n\ge 3$) be a domain in $(\R^n;u_1,u_2,\dots,u_n)$ and
\[
  f\colon{}U^n\longrightarrow (\R^{n+1},g_0)
\]
a front, that is, there exists a unit vector field $\nu$
(called the {\em unit normal vector})
such that $g_0(f_*X,\nu)=0$ for all $X\in TU^n$
and $(f,\nu)\colon{}U^n\to \R^{n+1}\times S^n$ is an immersion.
We set
\[
  \lambda:=\det(f_{u_1},\dots,f_{u_n},\nu),
\]
and call it the {\em signed volume density function}.
A point $p\in U^n$ is called a {\em singular point\/} if
$f$ is not an immersion at $p$.
Moreover, if $d\lambda\ne 0$ at $p$, we call $p$  a
{\em non-degenerate singular point}.
On a sufficiently small neighborhood of a non-degenerate singular point
$p$, 
the singular set is a $(n-1)$-dimensional submanifold called
the {\em singular submanifold}.
The $1$-dimensional vector space at the  non-degenerate singular point $p$
which is the kernel of the differential map 
$(f_*)_p\colon{}T_pU^n\to \R^{n+1}$ is called
the {\em null direction}. 
We call $p\in U^n$  a {\em cuspidal edge\/} 
if the null direction is transversal to the singular submanifold.
Then, by a similar argument to the proof of 
Fact~\ref{fact:intrinsic-criterion} in \cite{KRSUY},
one can prove that a cuspidal edge is an $A_2$-singularity, that is, 
locally diffeomorphic at the origin to the front
$f_C(u_1,\dots,u_n)=(u_1^2,u_1^3,u_2,\dots,u_n)$.
\begin{theorem}
 \label{thm:hyper}
 Let $f\colon{}U^n\to (\R^{n+1},g_0)$ $(n\ge 3)$ be a front whose 
 singular points are all cuspidal edges. 
 If the sectional curvature $K$ at the regular points
 is bounded, then the second fundamental form on the singular submanifold
 vanishes. 
 Moreover, if $K$ is positive everywhere on the regular set,
 the sectional curvature of the singular submanifold is non-negative.
 Furthermore, if $K\geq \delta(>0)$, then
 the sectional curvature of the singular submanifold is positive.
\end{theorem}
\begin{remark}
 The previous Theorem \ref{thm:gaussian-singular} is  
 deeper than this theorem.
When $n\ge 3$ we can 
consider sectional curvature on the singular set, but when 
$n=2$ the singular set is $1$-dimensional and so we cannot define
the sectional curvature.
Rather, one defines the singular curvature instead.
We do not define singular curvature for fronts
 when $n\ge 3$.
\end{remark}
\begin{proof}[Proof of Theorem~\ref{thm:hyper}]
 Without loss of generality, we may assume that
 the singular submanifold of $f$ is the $(u_1,\dots,u_{n-1})$-plane,
 and $\partial_n:=\partial/\partial {u_n}$ is the null direction.
 To prove the first assertion, it is sufficient to show that
 $h(X,X)=0$ for an arbitrary fixed tangent vector of the 
 singular submanifold.
 By changing coordinates if necessary, we may assume that
 $X=\partial_1=\partial/\partial u_1$.
 The sectional curvature $K(\partial_1\wedge \partial_n)$
 with respect to the $2$-plane spanned by 
 $\{\partial_1,\partial_{n}\}$ is given by
 \[
   K(\partial_1\wedge \partial_n)=
   \frac{h_{11}h_{nn}-(h_{1n})^2}{g_{11}g_{nn}-(g_{1n})^2}\qquad
   \bigl(
    g_{ij}=g_0(\partial_i,\partial_j),~
     h_{ij}=h(\partial_i,\partial_j)
   \bigr),
 \]
 where $h$ is the second fundamental form.
 By the same reasoning as in the proof of 
 Theorem~\ref{thm:gaussian-singular},
 the boundedness of $K(\partial_1\wedge \partial_n)$ implies
 \[
   0=\left. \left (h_{11}h_{nn}-(h_{1n})^2\right)_{u_n}
    \right |_{u_n=0}\!\!\!
   =h_{11}
    \left.\frac{\partial h_{nn}}{\partial u_n}\right |_{u_n=0}
    \!\!\!
   =h_{11}
    \left. g_0(D_{u_n}f_{u_n},\nu_{u_n})\right |_{u_n=0}.
 \]
 To show $h_{11}=h(X,X)=0$, it is sufficient to show
 $g_0(D_{u_n}f_{u_n},\nu_{u_n})$ does not vanish when  $u_n=0$.
 Since $f$ is a front with non-degenerate singularities, 
 we have
 \[
   0\ne \lambda_{u_n}=\det(f_{u_1},\dots,f_{u_{n-1}},D_{u_n}f_{u_n},\nu),
 \]
 which implies 
 $f_{u_1},\dots,f_{u_{n-1}},D_{u_n}f_{u_n}$, and $\nu$
 are linearly independent when $u_n=0$, and then 
 $\nu_{u_n}$ can be written as a linear combination of them. 
 Since $f$ is a front, $\nu_{u_n}\ne 0$ holds when $u_n=0$, 
 and we have $2g_0(\nu_{u_n},\nu)=g_0(\nu,\nu)_{u_n}=0$,  and
 \[
    g_0(\nu_{u_n},f_{u_j})=-g_0(\nu,D_{u_n}f_{u_j})=
    g_0(\nu_{u_j},f_{u_n})=0
    \qquad (j=1,\dots,n-1).
 \]
 Thus we have that
 $g_0(D_{u_n}f_{u_n},\nu_{u_n})$ never vanishes at $u_n=0$.

 Next we show the non-negativity of the sectional curvature $K_S$
 of the singular manifold.
 It is sufficient to show $K_S(\partial_1\wedge \partial_2)\ge 0$
 at $u_n=0$.
 Since the sectional curvature $K_{U_n}$ is non-negative, 
 we have 
 \begin{equation}
  \label{eq:h1}
  \left. \frac{\partial^2}{(\partial u_n)^2}
   (h_{11}h_{22}-(h_{12})^2)\right |_{u_n=0}\ge 0,
 \end{equation}
by the same argument as in the proof of 
 Theorem~\ref{thm:gaussian-singular}.
 Since the restriction of $f$ to the singular manifold is
 an immersion, the Gauss equation yields that
 \[
   K_S(\partial_1\wedge \partial_2)
     =\frac{g_0(\alpha_{11},\alpha_{22})-g_0(\alpha_{12},\alpha_{12})}
       {g_{11}g_{22}-(g_{12})^2},
 \]
 where $\alpha$ is the second fundamental form of the
 singular submanifold in $\R^{n+1}$ and
 $\alpha_{ij}=\alpha(f_{u_j},f_{u_j})$.

 On the other hand, since the second fundamental form
 $h$ of $f$ vanishes,
 $g_0(\nu_{u_n},f_{u_j})=0$ holds for $j=1,\dots,n$,
 that is, $\nu$ and $\nu_{u_n}$ are linearly independent 
 vectors.
 Moreover, we have
 \begin{align*}
   \alpha_{ij}
       &=g_0(\alpha_{ij},\nu)\nu+
        \frac{1}{|\nu_{u_n}|^2}g_0(\alpha_{ij},\nu_{u_n})\nu_{u_n}
        \\
        &=h_{ij}\nu+
             \frac{1}{|\nu_{u_n}|^2} g_0(\alpha_{ij},\nu_{u_n})\nu_{u_n}
        =\frac{1}{|\nu_{u_n}|^2}(h_{ij})_{u_n}\nu_{u_n},
 \end{align*}
 since the second fundamental form
 $h$ of $f$ vanishes and
 \begin{align*}
  g_0(\alpha_{ij},\nu_{u_n})=
  g_0(D_{u_j}f_{u_i},\nu_{u_n})
  =(h_{ij})_{u_n}-g_0(D_{u_i}D_{u_j}f_{u_n},\nu)=(h_{ij})_{u_n}
 \end{align*}
 for $i,j=1,\dots,n-1$.
 Thus we have
 \[
   K_S(\partial_1\wedge \partial_2)
    =\frac{1}
     {g_{11}g_{22}-(g_{12})^2}
    \left. \frac{\partial^2}{(\partial u_n)^2}
    (h_{11}h_{22}-(h_{12})^2)\right |_{u_n=0}\ge 0.
 \]
\end{proof}
\begin{example}\label{ex:hyper}
 We set
 \[
   f(u,v,w):=(v,w,u^2+a v^2+b w^2,u^3+c u^2):\R^3\to \R^4,
 \]
 which gives a front with the unit normal vector 
 \begin{multline*}
   \nu = \frac{1}{\delta}\bigl(2av(2c+3u),2bw(2c+3u),-2c-3u,2\bigr),\\
  \qquad\text{where}\quad
   \delta=\sqrt{4+(3u+2c)^2(1+4a^2v^2+4b^2w^2)}.
 \end{multline*}
 The singular set is the $vw$-plane and the $u$-direction is
 the null direction.
 Then all singular points are cuspidal edges.
 The second fundamental form is given by 
 $h = \delta^{-1}\{6u\,du^2-2(3u+2c)(a\,dv^2+b\,dw^2)\}$,
 which vanishes on the singular set if $ac=bc=0$.

 On the other hand, the sectional curvatures are computed as
 \begin{align*}
    K(\partial_u\wedge\partial_v) &=
    \frac{12a(3u+2c)}{u\delta^2\bigl(4+(3u+2c)^2(1+4a^2v^2)\bigr)},
    \\
    K(\partial_u\wedge\partial_w) &=
    \frac{12b(3u+2c)}{u\delta^2\bigl(4+(3u+2c)^2(1+4b^2w^2)\bigr)},
 \end{align*}  
 which are bounded in a neighborhood of  the singular set if
 and only if $ac=bc=0$.
 If $ac=bc=0$, $K\geq 0$ if and only if $a\geq 0$ and $b\geq 0$, which 
 implies 
 $K_S=4ab(3u+2c)^2/(\delta^2|\partial_v\wedge\partial_w|^2)>0$.
\end{example}

\section{Intrinsic formulation}
\label{sec:intrinsic}

The Gauss-Bonnet theorem is intrinsic in nature,
and it it quite natural to formulate the singularities 
of wave fronts intrinsically.
We can characterize the limiting tangent bundles of the fronts
and can give the following abstract definition:
\begin{definition}
 Let $M^2$ be a $2$-manifold.
 An orientable vector bundle $\E$ of rank $2$
 with a metric $\inner{~}{~}$ and a metric connection
 $D$ is called {\it an abstract limiting tangent bundle\/} 
 or {\em a coherent tangent bundle\/} if
 there is a bundle homomorphism
 \[
   \psi\colon{}TM^2\longrightarrow \E
 \]
 such that
 \begin{equation}
  D_{X}\psi(Y)-D_{Y}\psi(X)=\psi([X,Y]) \qquad (X,Y\in TM^2).
 \end{equation}
\end{definition}
In this setting, the pull-back of the metric
$ds^2:=\psi^*\inner{~}{~}$ is called 
{\em the first fundamental form\/} of $\E$.
A point $p\in M^2$ is called a {\it singular point\/} if
the first fundamental form is not positive definite.
Since $\E$ is orientable, there exists a skew-symmetric
bilinear form 
$\mu_p\colon{}\E_p\times\E_p\to \R$ for each $p\in M^2$,
where $\E_p$ is the fiber of $\E$ at $p$,
such that $\mu(e_1,e_2)=\pm 1$ for any orthonormal frame $\{e_1,e_2\}$
on $\E$.

A frame $\{e_1,e_2\}$ is called positive if
$\mu(e_1,e_2)=1$.
A singular point $p$ is called {\it non-degenerate} if
the derivative $d\lambda$ of
the function
\begin{equation}\label{eq:intlambda}
  \lambda:=\mu\left(
           \psi\left(\frac{\partial}{\partial u}\right),
           \psi\left(\frac{\partial}{\partial v}\right)
           \right)
\end{equation}
does not vanish at $p$, 
where $(U;u,v)$ is a local coordinate system of $M^2$ at $p$.
On a neighborhood of a non-degenerate singular point,
the singular set consists of a regular curve, called the 
{\it singular curve}.
The tangential direction of the singular curve is
called the {\it singular direction}, and the direction of the kernel of
$\psi$ is called the {\em  null direction}.
Then we can define {\it intrinsic cuspidal edges\/} and
{\it intrinsic swallowtails\/} according to 
Fact~\ref{fact:intrinsic-criterion}.
For a given singular curve $\gamma(t)$ consisting of 
intrinsic cuspidal edge points, the singular curvature function
is defined by
\[
  \kappa_s(t):=\sign\bigl(\lambda(\eta)\bigr)\hat \kappa_g(t),
\]
where 
$\hat\kappa_g(t):=\inner{D_t\psi(\gamma'(t))}{n(t)}$
is the {\it limiting geodesic curvature},
$n(t)\in \E_{\gamma(t)}$ is a unit vector such that 
$\mu\bigl(\psi(\gamma'(t)),n(t)\bigr)=1$,
and $\eta(t)$ is the null direction such that 
$\bigl(\gamma'(t),\eta(t)\bigr)$ is a positive frame on $M^2$.
Then Theorem~\ref{thm:invariance-singular-curvature}
and Proposition~\ref{prop:intrinsic-singular-curvature} hold.
Let $(U;e_1,e_2)$ be an orthonormal frame field of $\E$
such that $\mu(e_1,e_2)=1$.
Then there exists a unique $1$-form $\alpha$ on $U$ such that
\[
 D_Xe_1=-\alpha(X)e_2,
 \qquad 
 D_Xe_2=\alpha(X)e_1\qquad (X\in TM^2),
\]
which is called the {\it connection form}.
Moreover, the exterior derivative $d\alpha$ does not
depend on the choice of a positive frame $(U;e_1,e_2)$
and gives a (globally defined) $2$-form on $M^2$.
When $M^2$ is compact, the integration
\begin{equation}\label{eq:Euler}
 \chi_{\E}^{}:=\frac{1}{2\pi}\int_{M^2} d\alpha 
\end{equation}
is an integer called the {\it Euler number\/} of $\E$.
Let $(U;e_1,e_2)$ be a positive orthonormal frame field of $\E$
and $\gamma(s)$ a curve in $U(\subset M^2)$ such that 
$\inner{\psi(\gamma'(s))}{\psi(\gamma'(s))}=1$.
Let $\phi(s)$ be the angle of $\psi(\gamma'(s))$ from
$e_1(\gamma(s))$.
Then we have
\begin{equation}\label{eq:GB-key}
 \hat\kappa_g\,ds=d\phi-\alpha.
\end{equation}
Let $\Delta$ be a triangle with interior angles $A,B,C$.
In the interior of $\Delta$, we suppose that there are no singular points
and that $\psi^*d\alpha$ is compatible with respect to the orientation
of $M^2$.
We give an orientation to $\partial \Delta$ such that 
conormal vector points into the domain $\Delta$.
By using the same argument as in the classical proof 
of the Gauss-Bonnet Theorem, 
we get the formulas \eqref{eq:GB-signed} and
\eqref{eq:GB-unsigned} in the introduction intrinsically.
This intrinsic formulation is meaningful if we consider the following
examples:
\begin{example}[Cuspidal cross caps]
 A map $f\colon{}M^2\to \R^3$ is called a {\it frontal\/} 
 if there exists a unit normal vector field $\nu$ 
 such that $f_*X$ is perpendicular to $\nu$ for all
 $X\in TM^2$. 
 A frontal is a front if $(f,\nu)\colon{}M^2\to \R^3\times S^2$
 is an immersion.
 A cuspidal cross cap is a singular point locally diffeomorphic to
 the map $(u,v)\mapsto (u,v^2,uv^3)$
 and is a frontal but not a front.
 In \cite{FSUY}, a useful criterion for cuspidal cross caps are given.
 Though a cuspidal cross cap is not a cuspidal edge, 
 the limiting tangent bundle is well defined and the singular point is
 an intrinsic cuspidal edge.
 In particular, our Gauss-Bonnet formulas hold for a frontal
 that admits only cuspidal edges, swallowtails and
 cuspidal cross caps,
 and degenerate peaks like as for a double swallowtail.
\end{example}
\begin{example}[Singularities with higher codimensions]
 A smooth map $f\colon{}M^2\to \R^n$ 
 defined on a $2$-manifold $M^2$ into $\R^n$ ($n>3$) 
 is called {\it an admissible map\/} if there exists a map
 $\nu\colon{}M^2\to G_2(\R^n)$
 into the oriented $2$-plane Grassman manifold $G_2(\R^n)$,
 such that it coincides with the Gauss map of $f$ on regular points
 of $f$.
 For an admissible map, the limiting tangent bundle is canonically 
 defined and we can apply our intrinsic formulation to it.
\end{example}
A realization problem for abstract limiting tangent bundles  
is  investigated in \cite{SUY2}.
The realization of first fundamental forms with singularities
has been treated in \cite{K2}.


\end{document}